# A MAXIMUM LIKELIHOOD METHOD FOR THE INCIDENTAL PARAMETER PROBLEM

By Marcelo J. Moreira

*Columbia University and FGV/EPGE*

This paper uses the invariance principle to solve the incidental parameter problem of [*Econometrica* **16** (1948) 1–32]. We seek group actions that preserve the structural parameter and yield a maximal invariant in the parameter space with fixed dimension. M-estimation from the likelihood of the maximal invariant statistic yields the maximum invariant likelihood estimator (MILE). Consistency of MILE for cases in which the likelihood of the maximal invariant is the product of marginal likelihoods is straightforward. We illustrate this result with a stationary autoregressive model with fixed effects and an agent-specific monotonic transformation model.

Asymptotic properties of MILE, when the likelihood of the maximal invariant does not factorize, remain an open question. We are able to provide consistent, asymptotically normal and efficient results of MILE when invariance yields Wishart distributions. Two examples are an instrumental variable (IV) model and a dynamic panel data model with fixed effects.

**1. Introduction.** The maximum likelihood estimator (MLE) is a procedure commonly used to estimate a parameter in stochastic models. Under regularity conditions, the MLE is not only consistent but also asymptotic optimal (e.g., [26]). In the presence of incidental parameters, however, the MLE of structural parameters may not be consistent. This failure occurs because the dimension of incidental parameters increases with the sample size, affecting the ability of MLE to consistently estimate the structural parameters. This is the so-called incidental parameter problem after the seminal paper by [35].

This paper appeals to the invariance principle to solve the incidental parameter problem. We propose to find a group action that preserves the model

Received June 2008; revised December 2008.

[1]Supported by the NSF Grant SES-0819761.

and the structural parameter. This yields a maximal invariant statistic. Its distribution depends on the parameters only through the maximal invariant in the parameter space. Maximization of the invariant likelihood yields the maximum invariant likelihood estimator (MILE). Distinct group actions in general yield different estimators. We seek group actions whose maximal invariant in the parameter space has fixed dimension regardless of the sample size.

The use of invariance to eliminate nuisance parameters has a long history (e.g., [9]). However, the use of invariance to solve the incidental parameter problem is limited to only a few models (e.g., see [29] for estimate variance components using invariance to the mean). There has also been some discussion on identifiability by [28] for additional groups of transformations. However, asymptotic properties of MILE are hardly addressed in the literature. The difficulty in obtaining asymptotic results arises because the likelihood of the maximal invariant is often not the product of marginal likelihoods.

An important methodological question is whether the use of invariance yields consistency and optimality in models whose number of parameters increases with the sample size. As is customary in the literature, we illustrate these results with a series of examples.

To establish a context, Section 3 considers two groups of transformations whose use of invariance completely discards the incidental parameters. In both examples, the likelihood of the maximal invariant is the product of marginal likelihoods; consistency, asymptotic normality, and efficiency of MILE are straightforward. The first example is the stationary autoregressive model with fixed effects. For a particular group action, the solution coincides with [4] conditional and [15] and [25] integrated likelihood approaches. The second example is the monotonic transformation model. The proposed transformation is agent-specific and has infinite dimension. The conditional and integrated likelihood approaches do not seem to be applicable here. The invariance principle provides an estimator that is consistent and asymptotically normal under the assumption of normal errors.

We then proceed to the two main examples of the paper. For both examples, invariance arguments yield Wishart distributions. Standardization of the likelihoods yields consistency, asymptotic normality, and optimality results for MILE. Although our theoretical findings are somewhat specific to Wishart distributions, we hope that interesting general lessons can be learned from studying those particular likelihoods.

Section 4 considers an instrumental variable (IV) model with $N$ observations and $K$ instruments. For the orthogonal group of transformations, MILE coincides with the LIMLK estimator. The asymptotic theory for the invariant likelihood unifies theoretical findings for LIMLK under both the strong instruments (SIV) and many weak instruments (MWIV) asymptotics (e.g., [10, 22] and [31]). This framework parallels standard M-estimation in problems in which the number of parameters does not change with the sample size. In particular, we are able to (i) show consistency of the MLE in the IV setup even under MWIV asymptotics from the perspective of likelihood maximization; (ii) derive the asymptotic distribution of the MLE directly from the objective function under SIV and MWIV asymptotics; and (iii) provide an explanation for optimality of MLE within the class of regular invariant estimators.

Section 5 presents a simple dynamic panel data model with $N$ individuals and $T$ time periods. We propose to use MILE based on the orthogonal group of transformations. This estimator is novel in the dynamic panel data literature and presents a number of desirable properties. It is consistent, as long as $NT$ goes to infinity (regardless of the relative rate of $N$ and $T$) and asymptotically normal under (i) large $N$, fixed $T$; and (ii) large $N$, large $T$ asymptotics when the autoregressive parameter is smaller than one. We derive an efficiency bound for large $N$, fixed $T$ asymptotics when errors are normal; our bound coincides with [17] bound when $T \to \infty$. MILE reaches (i) our bound when $N$ is large and $T$ is fixed; and (ii) [17] bound when both $N$ and $T$ are large. The bias-corrected ordinary least squares (BCOLS) estimator (e.g., [17]) only reaches the second bound. As a result, it is shown that MILE asymptotically dominates the BCOLS estimator. Finally, [13] use invariance to show that the correlated random effects estimator has a minimax property. The fixed effects estimator MILE also has a minimax property for the group of transformations considered here.

Section 6 compares MILE with existing fixed-effects estimators for the dynamic panel data model.

Section 7 concludes. The Appendix provides proofs for our results.

**2. The maximum invariant likelihood estimator.** In this section, we revisit the basic concepts of invariance (e.g., [16]) and their use to eliminate nuisance parameters. Let $P_{\gamma,\eta}$ denote the distribution of the data set $Y \in \mathbf{Y}$ when the structural parameter is $\gamma \in \mathbf{\Gamma}$ and the incidental parameter is $\eta \in \mathbf{N}$: $\mathcal{L}(Y) = P_{\gamma,\eta} \in \mathbf{P}$.

We seek a group $\mathbf{G}$ and actions $\mathcal{A}_1(\cdot, Y)$ and $\mathcal{A}_2(\cdot, (\gamma, \eta))$ in the sample and parameter spaces that preserve the model $\mathbf{P}$:

$$\mathcal{L}(Y) = P_{\gamma,\eta} \quad \Rightarrow \quad \mathcal{L}(\mathcal{A}_1(g, Y)) = P_{\mathcal{A}_2(g,(\gamma,\eta))} \qquad \text{for any } P_{\gamma,\eta} \in \mathbf{P}.$$

We are interested in $\gamma$. This yields the following definition.

DEFINITION 2.1. Suppose that $\mathcal{A}_2 \colon \mathbf{G} \times \mathbf{\Gamma} \times \mathbf{N} \to \mathbf{\Gamma} \times \mathbf{N}$ induces an action $\mathcal{A}_3 \colon \mathbf{G} \times \mathbf{N} \to \mathbf{N}$ such that

$$\mathcal{A}_2(g, (\gamma, \eta)) = (\gamma, \mathcal{A}_3(g, \eta)).$$

Then the parameter $\gamma$ is said to be preserved. The incidental parameter space $\mathbf{N}$ is preserved if

$$\mathbf{N} = \{\eta \in \mathbf{N}; \eta = \mathcal{A}_3(g, \widetilde{\eta}) \text{ for some } \widetilde{\eta} \in \mathbf{N}\}.$$

Suppose that both $\gamma$ and $\mathbf{N}$ are preserved. We can then appeal to the *invariance principle* and focus on invariant statistics $\phi(Y)$ in which $\phi(\mathcal{A}_1(g, Y)) = \phi(Y)$ for every $Y \in \mathcal{Y}$ and $g \in \mathbf{G}$. Any invariant statistic can be written as a function of a maximal invariant statistic defined below.

DEFINITION 2.2. A statistic $M \equiv M(Y)$ is a maximal invariant in the sample space if

$$M(\widetilde{Y}) = M(Y) \quad \text{if and only if} \quad \widetilde{Y} = \mathcal{A}_1(g, Y) \qquad \text{for some } g \in G.$$

An orbit of $\mathbf{G}$ is an equivalence class of elements $Y$, where $\widetilde{Y} \sim Y \pmod{\mathbf{G}}$, if there exists $g \in \mathbf{G}$ such that $\widetilde{Y} = \mathcal{A}_1(g, Y)$. By definition, $M$ is a *maximal invariant* statistic if it is invariant and takes distinct values on different *orbits* of $\mathbf{G}$. Every invariant procedure can be written as a function of a maximal invariant. Hence, we restrict our attention to the class of decision rules that depend only on the maximal invariant statistic. An analogous definition holds for the parameter space.

DEFINITION 2.3. A parameter $\theta \equiv \theta(\gamma, \eta)$ is a maximal invariant in the parameter space if $\theta(\gamma, \eta)$ is invariant and takes different values on different *orbits* of $\mathbf{G}$: $O_{\gamma,\eta} = \{\mathcal{A}_2(g, (\gamma, \eta)) \in \mathbf{\Gamma} \times \mathbf{N}; \text{for some } g \in \mathbf{G}\}$.

The distribution of a maximal invariant $M$ depends on $(\gamma, \eta)$ only through $\theta$. If $\mathcal{A}_2 : \mathbf{G} \times \mathbf{\Gamma} \times \mathbf{N} \to \mathbf{\Gamma} \times \mathbf{N}$ induces a group action $\mathcal{A}_3 : \mathbf{G} \times \mathbf{N} \to \mathbf{N}$, then $\theta \equiv (\gamma, \lambda)$, where $\lambda \in \mathbf{\Lambda}$ is the maximal invariant in the nuisance parameter space $\mathbf{N}$. The parameter set $\mathbf{\Lambda}$ is allowed to be the empty set.

DEFINITION 2.4. Let $f(M; \theta)$ be the p.d.f./p.m.f. of a maximal invariant statistic (we shall abbreviate $f(M; \theta)$ as the invariant likelihood). The maximum invariant likelihood estimator (MILE) is defined as

$$\widehat{\theta} \equiv \arg\max_{\widetilde{\theta} \in \Theta} f(M; \theta).$$

*Comments.* 1. Hereinafter, we assume the set $\Theta$ to be compact.

2. The estimator $\widehat{\theta}$ is the same for any one-to-one transformation of $M$. Different group actions $\mathcal{A}_1(\cdot, Y)$ and $\mathcal{A}_2(\cdot, (\gamma, \eta))$, however, yield different estimators. Hence, a better notation for $\widehat{\theta}$ would indicate its dependence on the choice of group actions.

3. In general, we seek group actions $\mathcal{A}_1(\cdot, Y)$ and $\mathcal{A}_2(\cdot, (\gamma, \eta))$ that preserve the model $\mathbf{P}$ and the structural parameter $\gamma$, and yield a maximal invariant $\lambda$ in $\mathbf{N}$ which has fixed dimension with the sample size.

We introduce some additional notation. The superscript $*$ indicates the true value of a parameter (e.g., $\gamma^*$ is the true value of the structural parameter $\gamma$). The subscript $N$ denotes dependence on the sample size $N$ (e.g., $\lambda_N^*$ is the true value of the maximal invariant $\lambda$ when the sample size is $N$). In addition, let $1_T$ be a $T$-dimensional vector of ones, $O_{j \times k}$ be a $j \times k$ matrix with entries zero, $e_j$ be a vector with entry $j$ equals one and other entries zero.

Hereinafter, additional notation is specific to each example.

**3. Transformations within individuals.** In this section, we present two examples of transformations within individuals. Instead of $P_{\gamma, \eta}$, we work with $P_{\gamma, \eta_i}^i$, the probability of the model for agent $i$. This clarifies our exposition and highlights the fact that the likelihood of each maximal invariant $M = (M_1, \ldots, M_N)$ is the sum of marginal likelihoods. In all examples below, the maximal invariant in the parameter space is $\theta = \gamma$, with the objective function simplifying to

$$(3.1) \qquad Q_N(\theta) = \frac{1}{N} \sum_{i=1}^{N} \ln f_i(m_i; \theta),$$

where $f_i(m_i; \theta)$ is the marginal density of the maximal invariant $M_i$ for each individual $i$. Because the MILE $\widehat{\theta}_N$ maximizes $Q_N(\theta)$, consistency, asymptotic normality and optimality of $\widehat{\theta}_N$ follow from standard results.

LEMMA 3.1. *Let $Q_N(\theta)$ be defined as in (3.1) and take all limits as $N \to \infty$.*

(a) *Suppose that* (i) $\sup_{\theta \in \Theta} |Q_N(\theta) - Q(\theta)| \to_p 0$ *for a fixed, nonstochastic function $Q(\theta)$, and* (ii) $\forall \varepsilon > 0$, $\inf_{\theta \notin B(\theta^*, \varepsilon)} Q(\theta) > Q(\theta^*)$. *Then*

$$\widehat{\theta}_N \to_p \theta^*.$$

(b) *Suppose that* (i) $\widehat{\theta}_N \to_p \theta^*$, (ii) $\theta^* \in \mathrm{int}(\Theta)$, (iii) $Q_N(\theta)$ *is twice continuously differentiable in some neighborhood of $\theta^*$,* (iv) $\sqrt{N} \partial Q_N(\theta^*)/\partial \theta \to_d N(0, \mathcal{I}(\theta^*))$, *and* (v) $\sup_{\theta \in \Theta} |\partial^2 Q_N(\theta^*)/\partial \theta \, \partial \theta' + \mathcal{I}(\theta)| \to_p 0$ *for some nonstochastic matrix that is continuous at $\theta^*$ where $\mathcal{I}(\theta^*)$ is nonsingular. Then*

$$\sqrt{N}(\widehat{\theta}_N - \theta^*) \to_d N(0, \mathcal{I}(\theta^*)^{-1}).$$

(c) *Suppose that* (i) $\{Q_N(\theta); \theta \in \Theta\}$ *is differentiable in quadratic mean at $\theta^*$ with nonsingular information matrix $\mathcal{I}(\theta^*)$, and* (ii) $\sqrt{N}(\widehat{\theta}_N - \theta^*) =$

$\mathcal{I}(\theta^*)^{-1}\sqrt{N}\partial Q_N(\theta^*)/\partial\theta + o_{Q_N(\theta^*)}(1)$. *Then*

$$\ln\frac{Q_N(\theta + h \cdot N^{-1/2})}{Q_N(\theta)} = h'S_N - \frac{1}{2}h'\mathcal{I}(\theta^*)h + o_{Q_N(\theta^*)}(1),$$

*where* $S_N \to_d N(0, \mathcal{I}(\theta^*))$ *under* $Q_N(\theta^*)$, *and* $\widehat{\theta}_N$ *is the best regular invariant estimator of* $\theta^*$.

*Comment.* Part (a) assumes (i) uniform convergence of $Q_N(\theta)$ and (ii) unique identifiability of $\theta^*$. Under the assumption that $\Theta$ is compact, [7] show that $Q_N(\theta) \to_p Q(\theta)$ uniformly, if and only if $Q_N(\theta) \to_p Q(\theta)$ pointwise, and $Q_N(\theta) - Q(\theta)$ is stochastically equicontinuous. The nonstochastic function $Q(\theta)$ satisfies the unique identifiability condition if $\theta$ is identified and $Q(\theta)$ is continuous.

**3.1. *A linear stationary panel data model.*** As an introductory example, consider a linear stationary panel data model with exogenous regressors and fixed effects:

$$y_{it} = \eta_i + x'_{it}\beta + u_{it},$$

where $y_{it} \in \mathbb{R}$ and $x_{it} \in \mathbb{R}^K$ are observable variables; $u_{it}$ are unobservable (possibly autocorrelated) errors, $i = 1, \ldots, N$, $t = 1, \ldots, T$; $\beta \in \mathbb{R}^K$ and $\sigma^2 \in \mathbb{R}$ are the structural parameters; and $\eta_i \in \mathbb{R}$ are incidental parameters, $i = 1, \ldots, N$.

The model for $y_{i\cdot} = [y_{i1}, \ldots, y_{iT}]' \in \mathbb{R}^T$ conditional on $x_{i\cdot} = [x_{i1}, \ldots, x_{iT}]' \in \mathbb{R}^{T \times K}$ is

$$y_{i\cdot} \overset{\text{ind}}{\sim} N(\eta_i 1_T + x_{i\cdot}\beta, \sigma^2 \Sigma_T)$$

(3.2)

$$\text{where } \Sigma_T = \frac{1}{1 - \rho^2}\begin{bmatrix} 1 & \rho & \cdots & \rho^{T-1} \\ \rho & 1 & & \\ \vdots & & \ddots & \\ \rho^{T-1} & & & 1 \end{bmatrix}.$$

Both the model and the structural parameter $\gamma = (\beta, \sigma^2, \rho)$ are preserved by translations $g \cdot 1_T$ (where $g$ is a scalar),

$$y_{i\cdot} + g \cdot 1_T \overset{\text{ind}}{\sim} N((\eta_i + g)1_T + x_{i\cdot}\beta, \sigma^2 \Sigma_T).$$

PROPOSITION 3.1. *Let* $g$ *be elements of the real line with* $g_1 \circ g_2 = g_1 + g_2$. *If the actions on the sample and parameter spaces are, respectively,* $\mathcal{A}_1(g, y_{i\cdot}) = (y_{i\cdot} + g \cdot 1_T)$ *and* $\mathcal{A}_2(g, (\beta, \sigma^2, \rho, \eta_i)) = (\beta, \sigma^2, \rho, \eta_i + g)$, *then:*

(a) *the vector $M_i = Dy_{i\cdot}$ is a maximal invariant in the sample space, where $D$ is a $T-1 \times T$ differencing matrix with typical row $(0, \ldots, 0, 1, -1, 0, \ldots, 0)$,*

(b) *$\gamma$ is a maximal invariant in the parameter space, and*

(c) *$M_i \overset{\text{ind}}{\sim} N(Dx_{i\cdot}\beta, \sigma^2 D\Sigma_T D')$ with density at $m_i = Dy_{i\cdot}$ given by*

$$f_i(m_i; \beta, \rho, \sigma^2) = (2\pi\sigma^2)^{-(T-1)/2} |D\Sigma_T D'|^{-1/2}$$
$$\times \exp\left\{ -\frac{1}{2\sigma^2}(y_{i\cdot} - x_{i\cdot}\beta)'D'(D\Sigma_T D')^{-1}D(y_{i\cdot} - x_{i\cdot}\beta) \right\}.$$

*Comment.* Under regularity conditions (e.g., (i) $\frac{1}{N}\sum_{i=1}^{N} \text{vec}(x_{i\cdot})\text{vec}(x_{i\cdot})' \to_p \Omega_{XX}$ p.d., (ii) $\frac{1}{\sqrt{N}}\sum_{i=1}^{N} u_{i\cdot} \otimes \text{vec}(x_{i\cdot}) \to_d N(0, \sigma^{*2}\Sigma_T^* \otimes \Omega_{XX})$, where $u_{i\cdot} = [u_{i1}, \ldots, u_{iT}]'$, (iii) $\sup_{N \geq 1} \frac{1}{N}\sum_{i=1}^{N} E\,\text{vec}(x_{i\cdot})\text{vec}(x_{i\cdot})' < \infty$, (iv) $(\beta, 1, 0) \notin \Theta$, $\forall \beta$, and (v) $\theta^* \in \text{int}(\Theta)$), we can use Lemma 3.1 to show that $\widehat{\theta}_N$ is consistent and asymptotically normal.

3.2. *A linear transformation model.* Consider a simple panel data transformation model,

$$\eta_i(y_{it}) = x_{it}'\beta + u_{it},$$

where $y_{it} \in \mathbb{R}$ and $x_{it} \in \mathbb{R}^K$ are observable variables; $u_{it} \in \mathbb{R}$ are unobservable errors, $i = 1, \ldots, N$, $t = 1, \ldots, T$, with $T > K$; $\eta_i : \mathbb{R} \to \mathbb{R}$ is an unknown, continuous, strictly increasing incidental function; and $\beta \in \mathbb{R}^K$ is the structural parameter. Unlike [2], we shall parameterize the distribution of the errors, $u_{it} \overset{\text{i.i.d.}}{\sim} N(\alpha_i, \sigma_i^2)$. Because of location and scale normalizations, we shall assume without loss of generality that $u_{it} \overset{\text{i.i.d.}}{\sim} N(0, 1)$.

The model for $y_{i\cdot} = (y_{i1}, y_{i2}, \ldots, y_{iT}) \in \mathbb{R}^T$ is then given by

$$P(y_{i\cdot} \leq v) = \prod_{t=1}^{T} \Phi(\eta_i(v_t) - x_{it}'\beta) \qquad \text{where } v = [v_1, v_2, \ldots, v_T]'.$$

Both the model and the structural parameter $\gamma \equiv \beta$ are preserved by continuous, strictly increasing transformations.

PROPOSITION 3.2. *Let $g$ be elements of the group of continuous, strictly increasing transformations, with $g_1 \circ g_2 = g_1(g_2)$. If the actions on the sample and parameter spaces are, respectively, $\mathcal{A}_1(g, (y_{i1}, y_{i2}, \ldots, y_{iT})) = (g(y_{i1}), g(y_{i2}), \ldots, g(y_{iT}))$ and $\mathcal{A}_2(g, (\beta, \eta_i)) = (\beta, \eta_i(g^{-1}))$, then:*

(a) *the statistic $M_i = (M_{i1}, \ldots, M_{iT})$ is the maximal invariant in the sample space, where $M_{it}$ is the rank of $y_{it}$ in the collection $y_{i1}, \ldots, y_{iT}$,*

(b) *the vector $\beta$ is the maximal invariant in the parameter space, and*

(c) $M_i$, $i = 1, \ldots, N$, are independent with marginal probability mass function of $M_i$ at $m_i$ given by

$$f_i(m_{i1}, \ldots, m_{iT}; \beta) = \frac{1}{T!} E\left[\exp\left\{\left(\sum_{t=1}^{T} V_{(m_{it})} x_{it}'\right)\beta\right\}\right]$$

$$\times \exp\left\{-\frac{1}{2}\beta'\left(\sum_{t=1}^{T} x_{it} x_{it}'\right)\beta\right\},$$

where $V_{(1)}, \ldots, V_{(T)}$ is an ordered sample from an $N(0, 1)$ distribution.

The likelihood of the maximal invariant also yields semiparametric methods. For example, consider the case in which $T = 2$. If $x_{i2}'\beta > x_{i1}'\beta$, then it is likely that $y_{i2} > y_{i1}$. This yields the semiparametric estimator of [2]. This estimator maximizes

$$Q_N(\beta) = \frac{1}{N} \sum_{i=1}^{N} \{H(y_{i2}, y_{i1})I(\triangle x_i'\beta > 0) + H(y_{i1}, y_{i2})I(\triangle x_i'\beta < 0)\},$$

where $H$ is an arbitrary function increasing in the first and decreasing in the second argument. This estimator is very appealing as it is consistent under more general error distributions. For asymptotic normality, [2] proposes to smoothen the objective function to obtain asymptotic normality whose convergence rate can be made arbitrarily close to $N^{-1/2}$. In contrast, the MILE estimator suggested here does not require arbitrary choices of $H$ or smoothening.

**4. An instrumental variables model.** Consider a simple simultaneous equations model with two endogenous variables, multiple instrumental variables (IVs) and errors that are normal with known covariance matrix. The model consists of a structural equation and a reduced-form equation:

$$y_1 = y_2\beta + u,$$

$$y_2 = Z\pi + v_2,$$

where $y_1, y_2 \in R^N$ and $Z \in R^{N \times K}$ are observed variables; $u, v_2 \in R^N$ are unobserved errors; and $\beta \in R$ and $\pi \in R^K$ are unknown parameters. The matrix $Z$ has full column rank $K$; the $N \times 2$ matrix of errors $[u : v_2]$ is assumed to be i.i.d. across rows with each row having a mean zero bivariate normal distribution with a nonsingular covariance matrix; $\pi$ is the incidental parameter; and $\beta$ is the parameter of interest.

The two-equation reduced-form model can be written in matrix notation as

$$Y = Z\pi a' + V,$$

where $Y = [y_1 : y_2]$, $V = [v_1 : v_2]$ and $a = (\beta, 1)'$. The distribution of $Y \in R^{N \times 2}$ is multivariate normal with mean matrix $Z\pi a'$, independence across rows and covariance matrix $\Sigma$ for each row.

Because the multivariate normal is a member of the exponential family of distributions, low-dimensional sufficient statistics are available for the parameter $(\beta, \pi')'$. Andrews, Moreira and Stock [8] and Chamberlain [12] propose using orthogonal transformations applied to the sufficient statistic $(Z'Z)^{-1/2}Z'Y$. The maximal invariant is $Y'N_Z Y$, where $N_Z = Z(Z'Z)^{-1}Z'$.

We shall use an invariance argument without reducing the data to a sufficient statistic. For convenience, it is useful to write the model in a canonical form. The matrix $Z$ has the polar decomposition $Z = \omega(\rho', 0_{K \times (N-K)})'$, where $\omega$ is an $N \times N$ orthogonal matrix, and $\rho$ is the unique symmetric, positive definite square root of $Z'Z$. Define $R = \omega'Y$ and let $\eta = \rho\pi$. Then the canonical model is

$$R \overset{d}{=} \begin{pmatrix} \eta a' \\ 0 \end{pmatrix} + V, \qquad \mathcal{L}(V) = N(0, I_N \otimes \Sigma).$$

Both model and structural parameters $\beta$ and $\Sigma$ are preserved by transformations $O(K)$ in the first $K$ rows of $R$. The next proposition obtains the maximal invariants in the sample and parameter spaces.

PROPOSITION 4.1. *Let $g$ be elements of the orthogonal group of transformations $O(K)$ and partition the sample space $R = (R_1', R_2')'$, where $R_1$ is $K \times 2$ and $R_2$ is $(N-K) \times 2$. If the actions on the sample and parameter spaces are, respectively, $\mathcal{A}_1(g, R) = ((gR_1)', R_2')'$ and $\mathcal{A}_2(g, (\beta, \Sigma, \eta)) = (\beta, \Sigma, g\eta)$, then:*

(a) *the maximal invariant in the sample space is $M = (R_1'R_1, R_2)$, and*

(b) *the maximal invariant in the parameter space is $\theta_N = (\beta, \Sigma, \lambda_N)$, where $\lambda_N \equiv \eta'\eta/N$.*

To illustrate the approach, we assume for simplicity that $\Sigma$ is known. Hence, we omit $\Sigma$ from now on [e.g., $\theta_N = (\beta, \lambda_N)$].

The density of $M$ is the product of the marginal densities of $R_1'R_1$ and $R_2$. Since $R_2$ is an ancillary statistic, we can focus on the marginal density of $R_1'R_1 \equiv Y'N_Z Y$ in the maximization of the log-likelihood. As the density of $Y'N_Z Y$ is not well-behaved as $N$ goes to infinity, we work with the density of $W_N \equiv N^{-1}Y'N_Z Y$ instead.

THEOREM 4.1. *The density of $W_N \equiv N^{-1}Y'N_Z Y$ evaluated at $w$ is*

$$g(w; \beta, \lambda_N) = C_{1,K} \cdot N^K \cdot \exp\left(-\frac{N\lambda_N}{2}a'\Sigma^{-1}a\right)|\Sigma|^{-K/2}|w|^{(K-3)/2}$$

$$\times \exp\left(-\frac{N}{2}\operatorname{tr}(\Sigma^{-1}w)\right)$$

(4.1)
$$\times (N\sqrt{\lambda_N \cdot a'\Sigma^{-1}w\Sigma^{-1}a})^{-(K-2)/2}$$
$$\times I_{(K-2)/2}(N\sqrt{\lambda_N \cdot a'\Sigma^{-1}w\Sigma^{-1}a}),$$

where $C_{1,K}^{-1} = 2^{(K+2)/2}\pi^{1/2}\Gamma(\frac{K-1}{2})$, $I_\nu(\cdot)$ denotes the modified Bessel function of the first kind of order $\nu$, and $\Gamma(\cdot)$ is the gamma function.

Define MILE as
$$\widehat{\theta}_N \equiv \arg\max_{\theta \in \Theta} Q_N(\theta),$$

where $Q_N(\theta) \equiv N^{-1}\ln g(W_N; \theta_N)$ and $\theta_N = (\beta, \lambda_N)$.[1] The next result shows that $\widehat{\theta}_N = \theta_N^* + o_p(1)$ under general conditions.

THEOREM 4.2. (a) Under the assumption that $N \to \infty$ with $K$ fixed or $K/N \to 0$, (i) if $\lambda_N^*$ is fixed at $\lambda^* > 0$, then $\widehat{\theta}_N \to_p \theta^* = (\beta^*, \lambda^*)$, (ii) if $\lambda_N^* \to_p \lambda^* > 0$, then $\widehat{\theta}_N \to_p \theta^* = (\beta^*, \lambda^*)$ and (iii) if $0 < \liminf \lambda_N^* \le \limsup \lambda_N^* < \infty$, then $\widehat{\theta}_N = \theta_N^* + o_p(1)$.

(b) Under the assumption that $N \to \infty$ with $K/N \to \alpha > 0$, (i) if $\lambda_N^*$ is fixed at $\lambda^* > 0$, then $\widehat{\theta}_N \to_p \theta^* = (\beta^*, \lambda^*)$, (ii) if $\lambda_N^* \to_p \lambda^* > 0$, then $\widehat{\theta}_N \to_p \theta^* = (\beta^*, \lambda^*)$ and (iii) if $0 < \liminf \lambda_N^* \le \limsup \lambda_N^* < \infty$, then $\widehat{\theta}_N = \theta_N^* + o_p(1)$, where $\theta_N^* = (\beta^*, \lambda_N^*)$.

*Comments.* 1. Parts (a), (b)(i) yield consistency results conditional on $\lambda_N^*$; the remaining results of the theorem are unconditional on $\lambda_N^*$. Parts (a), (b)(ii) yield consistency results for $\beta^*$ under SIV and MWIV asymptotics when $\lambda_N^* \to_p \lambda^*$. The assumption of $\lambda_N^* \to_p \lambda^*$ is standard in the literature, but parts (a), (b)(iii) show that $\widehat{\beta}_N \to_p \beta_N^*$ without imposing convergence of $\lambda_N^*$.

2. This result also holds under nonnormal errors, as long as $V(W_N) \to 0$.

PROPOSITION 4.2. *MILE of $\beta$ is the limited information maximum likelihood (LIMLK) estimator.*

Proposition 4.2 together with Theorem 4.2 explain why the LIMLK estimator is consistent when the number of instruments increases. The MILE estimator maximizes a log-likelihood function that is well-behaved as it depends on a finite number of parameters. The LIMLK estimator is consistent because it coincides with MILE.

---

[1] The objective function $Q_N(\theta)$ is not defined if $W_N$ is not positive definite (due to the term $\ln|W_N|$). To avoid this technical issue, we can instead maximize only the terms of $Q_N(\theta)$ that depend on $\theta$.

Theorem 4.3. *Let the score statistic and the Hessian matrix be*

$$S_N(\theta) = \frac{\partial \ln Q_N(\theta)}{\partial \theta} \quad and \quad H_N(\theta) = \frac{\partial^2 \ln Q_N(\theta)}{\partial \theta \, \partial \theta'},$$

*respectively, and define the matrix*

$$\mathcal{I}_\alpha(\theta^*) = \left[ \begin{array}{c} \lambda^{*2} \dfrac{a^{*\prime}\Sigma^{-1}a^* \cdot e_1'\Sigma^{-1}e_1(\alpha + 2\lambda^* a^{*\prime}\Sigma^{-1}a^*) + \alpha(a^{*\prime}\Sigma^{-1}e_1)^2}{(\alpha + \lambda^* a^{*\prime}\Sigma^{-1}a^*)(\alpha + 2\lambda^* a^{*\prime}\Sigma^{-1}a^*)} \\[2ex] \lambda^* \dfrac{a^{*\prime}\Sigma^{-1}e_1 \cdot a^{*\prime}\Sigma^{-1}a^*}{\alpha + 2\lambda^* a^{*\prime}\Sigma^{-1}a^*} \end{array} \right.$$

$$\left. \begin{array}{c} \lambda^* \dfrac{a^{*\prime}\Sigma^{-1}e_1 \cdot a^{*\prime}\Sigma^{-1}a^*}{\alpha + 2\lambda^* a^{*\prime}\Sigma^{-1}a^*} \\[2ex] \dfrac{(a^{*\prime}\Sigma^{-1}a^*)^2}{2(\alpha + 2\lambda^* a^{*\prime}\Sigma^{-1}a^*)} \end{array} \right].$$

(a) *Suppose that $\lambda_N^*$ is fixed at $\lambda^* > 0$ and $N \to \infty$ with $K$ fixed. Then (i) $\sqrt{N}S_N(\theta^*) \to_d N(0, \mathcal{I}_0(\theta^*))$, (ii) $H_N(\theta^*) \to_p -\mathcal{I}_0(\theta^*)$, and (iii) $\sqrt{N}(\widehat{\theta}_N - \theta^*) \to_d N(0, \mathcal{I}_0(\theta^*)^{-1})$.*

(b) *Suppose that $\lambda_N^*$ is fixed at $\lambda^* > 0$ and $N \to \infty$ with $K/N \to \alpha$. Then (i) $\sqrt{N}S_N(\theta^*) \to_d N(0, \mathcal{I}_\alpha(\theta^*))$, (ii) $H_N(\theta^*) \to_p -\mathcal{I}_\alpha(\theta^*)$ and (iii) $\sqrt{N}(\widehat{\theta}_N - \theta^*) \to_d N(0, \mathcal{I}_\alpha(\theta^*)^{-1})$.*

*Comment.* For convenience, we provide asymptotic results only for the case in which $\lambda_N^*$ is fixed at $\lambda^* > 0$. Small changes in the proofs also yield asymptotic results for $\lambda_N^* \to_p \lambda^*$.

As a corollary, we find the limiting distribution of LIMLK. This result coincides with those obtained by [10].

Corollary 4.1. *Define $\sigma_u^2 = b'\Sigma b$. Under SIV asymptotics (or under MWIV asymptotics with $\alpha = 0$), conditional on $\lambda_N^* = \lambda^* > 0$,*

$$(4.2) \qquad\qquad \sqrt{N}(\widehat{\beta}_N - \beta^*) \to_d N\left(0, \frac{\sigma_u^2}{\lambda^*}\right).$$

*Under MWIV asymptotics, conditional on $\lambda_N^* = \lambda^* > 0$,*

$$(4.3) \qquad \sqrt{N}(\widehat{\beta}_N - \beta^*) \to_d N\left(0, \frac{\sigma_u^2}{\lambda^{*2}}\left\{\lambda^* + \alpha\frac{1}{a^{*\prime}\Sigma^{-1}a^*}\right\}\right).$$

*Comments.* 1. The limiting distribution given in (4.3) simplifies to the one given in (4.2) as $\alpha \to 0$.

2. Instead of using the invariant likelihood to obtain a minimum distance (MD) estimator, we could instead use only its first moment. Define

$$(4.4) \qquad \overline{m}(W_N; \theta_N) = \mathrm{vech}\left(\frac{R_1'R_1}{N}\right) - \mathrm{vech}\left(aa' \cdot \lambda_N + \frac{K}{N}\Sigma\right).$$

If $\lambda_N^* > 0$, then the following holds (for possibly nonnormal errors):

$$(4.5) \qquad E_{\theta_N^*}(\overline{m}(W_N; \theta)) = 0 \quad \text{if and only if} \quad \theta_N = \theta_N^*.$$

Because the number of moment conditions does not increase under SIV or MWIV asymptotics, we can show that the MD estimator based on (4.4) and (4.5) is consistent and asymptotically normal.

Finally, we obtain the following result under SIV and MWIV asymptotics in our setup.

THEOREM 4.4. *Define the log-likelihood ratio*

$$\Lambda_N(\theta^* + h \cdot N^{-1/2}, \theta^*) = N(Q_N(\theta^* + h \cdot N^{-1/2}) - Q_N(\theta^*)).$$

(a) *Under SIV asymptotics,*

$$(4.6) \quad \Lambda_N(\theta^* + h \cdot N^{-1/2}, \theta^*) = h' \sqrt{N} S_N(\theta^*) - \tfrac{1}{2} h' \mathcal{I}_0(\theta^*) h + o_{Q_N(\theta^*)}(1),$$

*where* $\sqrt{N} S_N(\theta^*) \to_d N(0, \mathcal{I}_0(\theta^*))$ *under* $Q_N(\theta^*)$.

(b) *Under MWIV asymptotics,*

$$(4.7) \quad \Lambda_N(\theta^* + h \cdot N^{-1/2}, \theta^*) = h' \sqrt{N} S_N(\theta^*) - \tfrac{1}{2} h' \mathcal{I}_\alpha(\theta^*) h + o_{Q_N(\theta^*)}(1),$$

*where* $\sqrt{N} S_N(\theta^*) \to_d N(0, \mathcal{I}_\alpha(\theta^*))$ *under* $Q_N(\theta^*)$.

*Furthermore, the LIMLK estimator is asymptotically efficient within the class of regular invariant estimators under both SIV and MWIV asymptotics.*

*Comments.* 1. The proof of [14] uses asymptotic results by [19] for Wishart distributions. The standard literature on limit of experiments instead typically provides expansions around the score (e.g., [27]). Theorem 4.3 shows that the score is asymptotically normal with variance given by the reciprocal of the inverse of the limit of the Hessian matrix. As the remainder terms are asymptotically negligible, (4.6) and (4.7) hold true.

2. Theorem 4.4 requires the assumption of normal errors. Anderson, Kunitomo and Matsushita [6] exploit the fact that $W_N$ involves double sums (in terms of $N$ and $K$) to obtain optimality results for nonnormal errors.

Under MWIV asymptotics, the LIMLK estimator achieves the bound $(\mathcal{I}_\alpha(\theta^*)^{-1})_{11}$. Under SIV asymptotics, the bound $(\mathcal{I}_0(\theta^*)^{-1})_{11}$ for regular invariant estimators of $\beta$ is the same as the one achieved by limit of experiments applied to the likelihood of $Y$. Hence, there is no loss of efficiency in focusing on the class of invariant procedures under SIV asymptotics.

**5. A nonstationary dynamic panel data model.** Consider a simple dynamic panel data model with fixed effects,

$$y_{i,t} = \rho y_{i,t-1} + \eta_i + u_{it},$$

where $y_{it} \in \mathbb{R}$ are observable variables and $u_{it} \overset{\text{i.i.d.}}{\sim} N(0, \sigma^2)$ are unobservable errors, $i = 1, \ldots, N$, $t = 1, \ldots, T$; $\eta_i \in \mathbb{R}$ are incidental parameters, $i = 1, \ldots, N$; $\gamma = (\rho, \sigma^2) \in \mathbb{R} \times \mathbb{R}$ are structural parameters; and $y_{i,0}$ are the initial values of the stochastic process. We seek inference conditional on the initial values $y_{i,0} = 0$.[2]

In its matrix form, we have

$$(5.1) \qquad [y_{\cdot 1}, y_{\cdot 2}, \ldots, y_{\cdot T}] = \rho [y_{\cdot 0}, y_{\cdot 1}, \ldots, y_{\cdot T-1}] + \eta 1_T' + [u_{\cdot 1}, u_{\cdot 2}, \ldots, u_{\cdot T}],$$

where $y_{\cdot t} = [y_{1,t}, y_{2,t}, \ldots, y_{N,t}]' \in \mathbb{R}^N$, $u_{\cdot t} = [u_{1,t}, u_{2,t}, \ldots, u_{N,t}]' \in \mathbb{R}^N$, and $\eta = [\eta_1, \ldots, \eta_N]' \in \mathbb{R}^N$. Solving (5.1) recursively yields

$$[y_{\cdot 1}, y_{\cdot 2}, \ldots, y_{\cdot T}] = \eta(B1_T)' + [u_{\cdot 1}, u_{\cdot 2}, \ldots, u_{\cdot T}]B'$$

$$(5.2) \qquad \text{where } B = \begin{bmatrix} 1 & & \\ \vdots & \ddots & \\ \rho^{T-1} & \cdots & 1 \end{bmatrix}.$$

The inverse of $B$ has a simple form,

$$B^{-1} \equiv D = I_T - \rho \cdot J_T, \qquad \text{where } J_T = \begin{bmatrix} 0_{T-1}' & 0 \\ I_{T-1} & 0_{T-1} \end{bmatrix}$$

and $0_{T-1}$ is a $(T-1)$-dimensional column vector with zero entries.

If individuals $i$ are treated equally, the coordinate system used to specify the vectors $y_{\cdot t}$ should not affect inference based on them. In consequence, it is reasonable to restrict attention to coordinate-free functions of $y_{\cdot t}$. Indeed, we find that orthogonal transformations preserve both the model given in (5.2) and the structural parameter $\gamma = (\rho, \sigma^2)$.

PROPOSITION 5.1. *Let $g$ be elements of the orthogonal group of transformations $O(N)$. If the actions on the sample and parameter spaces are, respectively, $\mathcal{A}_1(g, Y) = gY$ and $\mathcal{A}_2(g, (\rho, \sigma^2, \eta)) = (\rho, \sigma^2, g\eta)$, then:*

(a) *the maximal invariant in the sample space is $M = Y'Y$, and*

(b) *the maximal invariant in the parameter space is $\theta_N = (\gamma, \lambda_N)$, where $\lambda_N = \eta'\eta/(N\sigma^2)$.*

---

[2]We can assume that $y_{i,0} = 0$ by writing the model as

$$(y_{i,t} - y_{i,0}) = \rho(y_{i,t-1} - y_{i,0}) + (\eta_i - y_{i,0}(1-\rho)) + u_{it},$$

for example, [25].

*Comment.* If there is autocorrelation $\Sigma_T$ that is homogeneous across individuals, the maximal invariant $M$ remains the same. The covariance matrix, however, changes to $\Sigma = \sigma^2 B \Sigma_T B'$.

For convenience, we standardize the distribution of $M = Y'Y$.

THEOREM 5.1. *If $N \geq T$, the density of $W_N \equiv N^{-1} Y'Y$ at $w$ is*

$$
\begin{aligned}
g(w; \rho, \sigma^2, \lambda_N) = {}& C_{2,N} \cdot (\sigma^2)^{-NT/2} |w|^{(N-T-1)/2} \\
& \times \exp\left(-\frac{N}{2\sigma^2} \operatorname{tr}(DwD')\right) \exp\left(-\frac{NT}{2} \lambda_N\right) \\
& \times \left(N \sqrt{\lambda_N \frac{1'_T DwD' 1_T}{\sigma^2}}\right)^{-(N-2)/2} \\
& \times I_{(N-2)/2}\left(N \sqrt{\lambda_N \frac{1'_T DwD' 1_T}{\sigma^2}}\right) \cdot N^{NT/2},
\end{aligned}
$$

(5.3)

*where $C_{2,N}^{-1} = 2^{NT/2 - (N-2)/2} \pi^{T(T-1)/4} \prod_{i=1}^{T-1} \Gamma(\frac{N-i}{2})$.*

Define MILE as

$$\widehat{\theta}_N \equiv \arg\max_{\theta \in \Theta} Q_N(\theta),$$

where $Q_N(\theta) \equiv (NT)^{-1} \ln g(W_N; \rho, \sigma^2, \lambda)$ and $\theta_N = (\rho, \sigma^2, \lambda_N)$.[3] The next result shows that $\widehat{\theta}_N = \theta_N^* + o_p(1)$ under general conditions.

THEOREM 5.2. (a) *Under the assumption that $N \to \infty$ with $T$ fixed,* (i) *if $\lambda_N^*$ is fixed at $\lambda^*$, then $\widehat{\theta}_N \to_p \theta^* = (\rho^*, \sigma^{*2}, \lambda^*)$,* (ii) *if $\lambda_N^* \to_p \lambda^*$, then $\widehat{\theta}_N \to_p \theta^* = (\rho^*, \sigma^{*2}, \lambda^*)$ and* (iii) *if $\limsup \lambda_N^* < \infty$, then $\widehat{\theta}_N = \theta_N^* + o_p(1)$, where $\theta_N^* = (\rho^*, \sigma^{*2}, \lambda_N^*)$.*

(b) *Under the assumption that $T \to \infty$ and $|\rho^*| < 1$,* (i) *if $\lambda_N^*$ is fixed at $\lambda^*$, then $\widehat{\theta}_N \to_p \theta^* = (\rho^*, \sigma^{*2}, \lambda^*)$,* (ii) *if $\lambda_N^* \to_p \lambda^*$, then $\widehat{\theta}_N \to_p \theta^* = (\rho^*, \sigma^{*2}, \lambda^*)$ and* (iii) *if $\limsup \lambda_N^* < \infty$, then $\widehat{\theta}_N = \theta_N^* + o_p(1)$, where $\theta_N^* = (\rho^*, \sigma^{*2}, \lambda_N^*)$.*

*Comments.* 1. This result also holds under nonnormal errors.

2. This theorem implies that $\widehat{\rho}_N \to_p \rho^*$ under the assumption that $NT \to \infty$ (regardless of the growing rate of $N$ and $T$).

The next result derives the limiting distribution of MILE when $N \to \infty$.

---

[3]If $N < T$, $W_N$ is not absolutely continuous with respect to the Lebesgue measure. We will still maximize the pseudo-likelihood to find $\widehat{\theta}_N$.

THEOREM 5.3. *Suppose that $\sigma^{*2} > 0$ and $\lambda_N^*$ is fixed at $\lambda^* > 0$, and let the score statistic and the Hessian matrix be*

$$S_N(\theta) = \frac{\partial \ln Q_N(\theta)}{\partial \theta} \quad and \quad H_N(\theta) = \frac{\partial^2 \ln Q_N(\theta)}{\partial \theta \, \partial \theta'},$$

*respectively, and define the matrix*

$$\mathcal{I}_T(\theta^*) = \begin{bmatrix} h_{1,T} + h_{2,T} + h_{3,T} & \dfrac{\lambda^*}{2\sigma^{*2}} \dfrac{1_T' F 1_T}{T} & \dfrac{1 + \lambda^* T}{1 + 2\lambda^* T} \dfrac{1_T' F 1_T}{T} \\[2ex] \dfrac{\lambda^*}{2\sigma^{*2}} \dfrac{1_T' F 1_T}{T} & \dfrac{1}{2(\sigma^{*2})^2} + \dfrac{\lambda^*}{4\sigma^{*2}} \dfrac{2\lambda^* T}{1 + 2\lambda^* T} & \dfrac{1}{4\sigma^{*2}} \\[2ex] \dfrac{1 + \lambda^* T}{1 + 2\lambda^* T} \dfrac{1_T' F 1_T}{T} & \dfrac{1}{4\sigma^{*2}} & \dfrac{1}{4\lambda^*} \end{bmatrix},$$

*where $DB^* \equiv I_T + (\rho^* - \rho)F$ and the three terms in the $(1,1)$ entry of $\mathcal{I}_T(\theta^*)$ are*

$$h_{1,T} = \frac{\operatorname{tr}(FF')}{T} + \lambda^* \frac{1_T' F' F 1_T}{T}, \qquad h_{2,T} = \frac{2\sigma^{*2}}{(1 + 2\lambda^* T)} \frac{(1_T' F 1_T)^2}{T}$$

*and*

$$h_{3,T} = -\frac{\lambda^*}{1 + \lambda^* T} \left\{ \frac{1_T' F' F 1_T}{T} + \lambda^* \frac{(1_T' F 1_T)^2}{T} \right\}.$$

*As $N \to \infty$ with $T$ fixed,*

(a) (i) $\sqrt{NT} S_N(\theta) \to_d N(0, \mathcal{I}_T(\theta^*))$, (ii) $H_N(\theta^*) \to_p -\mathcal{I}_T(\theta^*)$ *and* (iii) $\sqrt{NT}(\widehat{\theta}_N - \theta^*) \to_d N(0, \mathcal{I}_T(\theta^*)^{-1})$, *and*

(b) *the log-likelihood ratio is*

$$
\begin{aligned}
(5.4) \qquad & \Lambda_N(\theta^* + h \cdot (NT)^{-1/2}, \theta^*) \\
& = NT(Q_N(\theta^* + h \cdot (NT)^{-1/2}) - Q_N(\theta^*)) \\
& = h' \sqrt{NT} S_N(\theta^*) - \tfrac{1}{2} h' \mathcal{I}_T(\theta^*) h + o_{Q_N(\theta^*)}(1),
\end{aligned}
$$

$\sqrt{NT} S_N(\theta^*) \to_d N(0, \mathcal{I}_T(\theta^*))$ *under $Q_N(\theta^*)$. Furthermore, $\widehat{\theta}_N$ is asymptotically efficient within the class of regular invariant estimators under large $N$, fixed $T$ asymptotics.*

*Comments.* 1. It is possible to extend parts (a)(i), (iii) to nonnormal errors by finding the appropriate asymptotic distribution of $\sqrt{NT} S_N(\theta^*)$.

2. The MILE estimator $\widehat{\rho}_N$ achieves the bound $(\mathcal{I}_T(\theta^*)^{-1})_{11}$ as $N \to \infty$, whereas the bias-corrected OLS estimator does not.

3. Instead of using the invariant likelihood to obtain an estimator, we could instead use only its first moment. Let $w_i = y_{i\cdot} y_{i\cdot}'$, where $y_{i\cdot} = [y_{i,1}, y_{i,2}, \ldots, y_{i,T}]' \in \mathbb{R}^T$, and define

$$(5.5) \qquad \overline{m}(W_N; \theta_N) = \operatorname{vech}(W_N - \sigma^2 \operatorname{vech}(B\{I_T + \lambda_N \cdot 1_T 1_T'\} B)).$$

Then the following holds:

$$(5.6) \qquad E_{\theta_N^*}(\overline{m}(W_N; \theta_N)) = 0 \quad \text{if and only if} \quad \theta_N = \theta_N^*.$$

In the IV model, the number of moment conditions does not increase with $N$ or $K$ (see comment 2 to Corollary 4.1). In the panel data model, the number $T(T+1)/2$ of moment conditions given in (5.6) increases (too quickly) with $T$. Therefore, consistency and semiparametric efficiency results (e.g., [3] and [34]) do not apply to (5.6) as $T \to \infty$. Instead, Hahn and Kuersteiner [17] cleverly use Hájek's convolution theorem to obtain an efficiency bound for normal errors as $T \to \infty$ for the stationary case $|\rho^*| < 1$. The bias-corrected OLS estimator of $\rho$ achieves [17] bound for large $N$, large $T$ asymptotics.

Our efficiency bound $(\mathcal{I}_T(\theta^*)^{-1})_{11}$ reduces to [17] bound when $T \to \infty$. This shows that there is no loss of efficiency in focusing on the class of invariant procedures under large $N$, large $T$ asymptotics.

COROLLARY 5.1.  *Under the assumption that $|\rho^*| < 1$, the efficiency bound given by the $(1,1)$ coordinate of the inverse of $\mathcal{I}_\infty(\theta^*)^{-1} \equiv (\lim_{T \to \infty} \mathcal{I}_T(\theta^*))^{-1}$ converges to [17] efficiency bound of $(1 - \rho^{*2})$ as $T \to \infty$.*

As a final result, the MILE estimator $\widehat{\rho}_N$ also achieves the bound $(\mathcal{I}_T(\theta^*)^{-1})_{11}$ for large $N$, large $T$ asymptotics.

THEOREM 5.4.  *Under the assumption that $N \geq T \to \infty$, $|\rho^*| < 1$, and $\lambda_N^*$ is fixed at $\lambda^* > 0$, (i) $\sqrt{NT} S_N(\theta) \to_d N(0, \mathcal{I}_\infty(\theta^*))$, (ii) $H_N(\theta^*) \to_p -\mathcal{I}_\infty(\theta^*)$ and (iii) $\sqrt{NT}(\widehat{\theta}_N - \theta^*) \to_d N(0, \mathcal{I}_\infty(\theta^*)^{-1})$.*

**6. Numerical results.**  This section illustrates the MILE approach for estimation of the autoregressive parameter $\rho$ in the dynamic panel data model described in Section 5. The numerical results are presented as means and mean squared errors (MSEs) based on 1000 Monte Carlo simulations. These results are also available for other fixed-effects estimators: Arellano–Bond (AB), Ahn–Schmidt (AS) and bias-corrected OLS (BCOLS) estimators.

We consider different combinations between short and large panels: $N = 5$, 10, 25, 100 and $T = 2, 3, 5, 10, 25, 100$.

Table 1 presents the initial design from which several variations are drawn.[4] This design assumes that $\eta_i^* \overset{\text{i.i.d.}}{\sim} N(0,4)$ (random effects), $u_{it} \overset{\text{i.i.d.}}{\sim} N(0,1)$ (normal errors) and $\rho^* = 0.5$ (positive autocorrelation). The value $\sigma^*$ is fixed at one for all designs.

---

[4]The full set of results for $\rho$, $\sigma^2$, and $\lambda_N$ using different designs are available at http://www.columbia.edu/~mm3534/.

MILE seems to be correctly centered around 0.5. Even in a very short panel with $N = 5$ and $T = 2$, its bias of 0.0408 is quite small. As $N$ and/or $T$ increases, its mean approaches 0.5. For example, for $N = 5$ and $T = 25$, the bias is around 0.0129; for $N = 25$ and $T = 2$, the simulation mean is around 0.0040. These numerical results support the theoretical finding that MILE is consistent, as long as $NT$ goes to infinity (regardless of the relative rate of $N$ and $T$). The BCOLS estimator seems to have smaller bias than the AB and AS estimators for small $N$ and large $T$. The AB and AS estimators have large bias with small $N$ and $T$, but their performance improves with large $N$ and small $T$.

MILE also seems to have smaller MSE than the other estimators. The AS estimator outperforms the AB estimator in terms of MSE. The BCOLS estimator has smaller MSE than AS. The MSE of the BCOLS estimator, however, does not decrease if $N$ increases but $T$ is held constant. For $T \geq 25$, its performance is comparable to that of MILE. This provides numerical support for the theoretical finding that both MILE and BCOLS reach our large $N$, large $T$ bound.

TABLE 1
*Performance of estimators for the autoregressive parameter $\rho$ (random effects, normal errors, and $\rho = 0.50$)*

| $T$ | $N$ | Mean | | | | MSE | | | |
|---|---|---|---|---|---|---|---|---|---|
| | | MILE | BCOLS | AB | AS | MILE | BCOLS | AB | AS |
| 2 | 5 | 0.4592 | 0.9651 | * | * | 0.1552 | 0.4602 | * | * |
| 2 | 10 | 0.4859 | 0.9500 | * | * | 0.0631 | 0.3109 | * | * |
| 2 | 25 | 0.4960 | 0.9523 | * | * | 0.0246 | 0.2394 | * | * |
| 2 | 100 | 0.4974 | 0.9474 | * | * | 0.0054 | 0.2083 | * | * |
| 3 | 5 | 0.4431 | 0.7695 | −0.0578 | 0.8642 | 0.0631 | 0.1607 | 516.8489 | 0.3823 |
| 3 | 10 | 0.4789 | 0.7903 | 0.9766 | 0.8954 | 0.0280 | 0.1165 | 153.1105 | 0.2559 |
| 3 | 25 | 0.4908 | 0.8008 | 0.5705 | 0.9389 | 0.0115 | 0.1045 | 4.7087 | 0.2219 |
| 3 | 100 | 0.4979 | 0.8068 | 0.5372 | 0.9632 | 0.0024 | 0.0975 | 0.0724 | 0.2204 |
| 5 | 5 | 0.4626 | 0.6469 | 0.1980 | 0.6541 | 0.0231 | 0.0538 | 0.2323 | 0.0991 |
| 5 | 10 | 0.4802 | 0.6657 | 0.2386 | 0.7162 | 0.0116 | 0.0422 | 0.2145 | 0.0820 |
| 5 | 25 | 0.4935 | 0.6702 | 0.3768 | 0.7940 | 0.0044 | 0.0347 | 0.0869 | 0.1002 |
| 5 | 100 | 0.4991 | 0.6799 | 0.4650 | 0.8667 | 0.0010 | 0.0336 | 0.0136 | 0.1371 |
| 10 | 5 | 0.4731 | 0.5505 | 0.0385 | 0.3753 | 0.0122 | 0.0158 | 52.4500 | 0.0747 |
| 10 | 10 | 0.4861 | 0.5660 | 0.3249 | 0.4518 | 0.0049 | 0.0107 | 0.0489 | 0.0437 |
| 10 | 25 | 0.4937 | 0.5717 | 0.3977 | 0.5763 | 0.0021 | 0.0074 | 0.0211 | 0.0294 |
| 10 | 100 | 0.4993 | 0.5736 | 0.4625 | 0.7223 | 0.0005 | 0.0060 | 0.0058 | 0.0550 |
| 25 | 5 | 0.4871 | 0.5128 | ** | ** | 0.0048 | 0.0055 | ** | ** |
| 25 | 10 | 0.4930 | 0.5151 | ** | ** | 0.0025 | 0.0025 | ** | ** |
| 25 | 25 | 0.4966 | 0.5180 | ** | ** | 0.0010 | 0.0013 | ** | ** |
| 25 | 100 | 0.4997 | 0.5184 | ** | ** | 0.0002 | 0.0006 | ** | ** |
| 100 | 5 | 0.4941 | 0.5014 | ** | ** | 0.0014 | 0.0013 | ** | ** |
| 100 | 10 | 0.4978 | 0.5018 | ** | ** | 0.0007 | 0.0007 | ** | ** |
| 100 | 25 | 0.4990 | 0.5001 | ** | ** | 0.0003 | 0.0003 | ** | ** |
| 100 | 100 | 0.4997 | 0.5015 | ** | ** | 0.0001 | 0.0001 | ** | ** |

(*) The estimator is not available for $T = 2$.
(**) Computational cost is prohibitive for large $T$.

TABLE 2
*Performance of estimators for the autoregressive parameter ρ (nonconvergent effects, normal errors, and ρ = 0.50)*

| | | Mean | | | | MSE | | | |
|---|---|---|---|---|---|---|---|---|---|
| $T$ | $N$ | MILE | BCOLS | AB | AS | MILE | BCOLS | AB | AS |
| 2 | 5 | 0.4770 | 1.0835 | * | * | 0.0818 | 0.5044 | * | * |
| 2 | 10 | 0.4911 | 1.1389 | * | * | 0.0196 | 0.4442 | * | * |
| 2 | 25 | 0.4989 | 1.1994 | * | * | 0.0037 | 0.4959 | * | * |
| 2 | 100 | 0.5000 | 1.2352 | * | * | 0.0002 | 0.5410 | * | * |
| 3 | 5 | 0.4773 | 0.8349 | 0.2500 | 0.9455 | 0.0346 | 0.1603 | 384.7828 | 0.3733 |
| 3 | 10 | 0.4908 | 0.9110 | 0.5705 | 0.9203 | 0.0087 | 0.1818 | 0.5864 | 0.2215 |
| 3 | 25 | 0.4981 | 0.9636 | 0.5160 | 0.8997 | 0.0013 | 0.2173 | 0.0173 | 0.1719 |
| 3 | 100 | 0.4992 | 0.9904 | 0.5013 | 0.8231 | 0.0001 | 0.2406 | 0.0009 | 0.1049 |
| 5 | 5 | 0.4727 | 0.6997 | 0.2452 | 0.7159 | 0.0165 | 0.0603 | 0.1766 | 0.0873 |
| 5 | 10 | 0.4918 | 0.7415 | 0.4475 | 0.7635 | 0.0043 | 0.0640 | 0.0339 | 0.0795 |
| 5 | 25 | 0.4991 | 0.7755 | 0.4912 | 0.7902 | 0.0007 | 0.0768 | 0.0046 | 0.0861 |
| 5 | 100 | 0.4997 | 0.7936 | 0.4988 | 0.7854 | 0.0000 | 0.0863 | 0.0002 | 0.0816 |
| 10 | 5 | 0.4789 | 0.5798 | −0.9436 | 0.4278 | 0.0080 | 0.0151 | 1721.7952 | 0.0516 |
| 10 | 10 | 0.4908 | 0.6104 | 0.4005 | 0.5980 | 0.0024 | 0.0148 | 0.0197 | 0.0281 |
| 10 | 25 | 0.5027 | 0.6326 | 0.4806 | 0.7370 | 0.0014 | 0.0180 | 0.0022 | 0.0583 |
| 10 | 100 | 0.5000 | 0.6452 | 0.4988 | 0.7765 | 0.0000 | 0.0211 | 0.0001 | 0.0765 |
| 25 | 5 | 0.4884 | 0.5157 | ** | ** | 0.0040 | 0.0042 | ** | ** |
| 25 | 10 | 0.4949 | 0.5330 | ** | ** | 0.0014 | 0.0027 | ** | ** |
| 25 | 25 | 0.4995 | 0.5464 | ** | ** | 0.0003 | 0.0024 | ** | ** |
| 25 | 100 | 0.4999 | 0.5562 | ** | ** | 0.0000 | 0.0032 | ** | ** |
| 100 | 5 | 0.4964 | 0.4994 | ** | ** | 0.0013 | 0.0014 | ** | ** |
| 100 | 10 | 0.4987 | 0.5038 | ** | ** | 0.0006 | 0.0005 | ** | ** |
| 100 | 25 | 0.4994 | 0.5076 | ** | ** | 0.0002 | 0.0002 | ** | ** |
| 100 | 100 | 0.5001 | 0.5119 | ** | ** | 0.0000 | 0.0002 | ** | ** |

(*) The estimator is not available for $T = 2$.
(**) Computational cost is prohibitive for large $T$.

Table 2 reports results for $\lambda_N^* = N$ (nonconvergent effects), normal errors and $\rho^* = 0.5$. Table 3 presents results for random effects, $u_{it} \overset{\text{i.i.d.}}{\sim} (\chi^2(1) - 1)/\sqrt{2}$ (nonnormal errors) and $\rho^* = 0.5$. In both cases, MILE continues to have smaller bias and MSE than the other estimators. This result is surprising with nonnormal errors as the AB and AS estimators could potentially dominate MILE when $N$ is large and $T$ is small.

Tables 4 and 5 differ from Table 1 only in the autoregressive parameter; respectively, $\rho^* = -0.5$ (negative autocorrelation) and $\rho^* = 1.0$ (integrated

TABLE 3

*Performance of estimators for the autoregressive parameter ρ (random effects, nonnormal errors, and ρ = 0.50)*

| | | Mean | | | | MSE | | | |
|---|---|---|---|---|---|---|---|---|---|
| *T* | *N* | MILE | BCOLS | AB | AS | MILE | BCOLS | AB | AS |
| 2 | 5 | 0.4520 | 0.9797 | * | * | 0.1430 | 0.5085 | * | * |
| 2 | 10 | 0.5024 | 0.9975 | * | * | 0.0869 | 0.3687 | * | * |
| 2 | 25 | 0.4993 | 0.9665 | * | * | 0.0414 | 0.2711 | * | * |
| 2 | 100 | 0.5042 | 0.9507 | * | * | 0.0105 | 0.2175 | * | * |
| 3 | 5 | 0.4666 | 0.7910 | 0.3562 | 0.8923 | 0.0687 | 0.1811 | 31.5729 | 0.4008 |
| 3 | 10 | 0.4803 | 0.8056 | 0.4189 | 0.9204 | 0.0343 | 0.1373 | 59.3092 | 0.2723 |
| 3 | 25 | 0.4951 | 0.8054 | 0.3363 | 0.9376 | 0.0143 | 0.1104 | 53.3848 | 0.2233 |
| 3 | 100 | 0.4992 | 0.8091 | 0.5244 | 0.9683 | 0.0030 | 0.0999 | 0.0839 | 0.2278 |
| 5 | 5 | 0.4712 | 0.6629 | 0.2628 | 0.6585 | 0.0268 | 0.0647 | 0.1905 | 0.1359 |
| 5 | 10 | 0.4821 | 0.6704 | 0.3211 | 0.6975 | 0.0150 | 0.0456 | 0.1282 | 0.0872 |
| 5 | 25 | 0.4928 | 0.6778 | 0.3899 | 0.7748 | 0.0045 | 0.0380 | 0.0810 | 0.0914 |
| 5 | 100 | 0.4967 | 0.6798 | 0.4717 | 0.8539 | 0.0011 | 0.0339 | 0.0128 | 0.1291 |
| 10 | 5 | 0.4722 | 0.5602 | 0.0781 | 0.3906 | 0.0110 | 0.0175 | 162.8453 | 0.0840 |
| 10 | 10 | 0.4893 | 0.5663 | 0.3471 | 0.4507 | 0.0047 | 0.0105 | 0.0405 | 0.0516 |
| 10 | 25 | 0.4946 | 0.5721 | 0.4084 | 0.5625 | 0.0020 | 0.0077 | 0.0178 | 0.0309 |
| 10 | 100 | 0.4984 | 0.5745 | 0.4740 | 0.7154 | 0.0005 | 0.0061 | 0.0035 | 0.0514 |
| 25 | 5 | 0.4819 | 0.5113 | ** | ** | 0.0052 | 0.0046 | ** | ** |
| 25 | 10 | 0.4890 | 0.5157 | ** | ** | 0.0024 | 0.0026 | ** | ** |
| 25 | 25 | 0.4974 | 0.5182 | ** | ** | 0.0010 | 0.0014 | ** | ** |
| 25 | 100 | 0.4990 | 0.5187 | ** | ** | 0.0003 | 0.0006 | ** | ** |
| 100 | 5 | 0.4949 | 0.4997 | ** | ** | 0.0015 | 0.0014 | ** | ** |
| 100 | 10 | 0.4972 | 0.5004 | ** | ** | 0.0007 | 0.0007 | ** | ** |
| 100 | 25 | 0.5000 | 0.5015 | ** | ** | 0.0003 | 0.0003 | ** | ** |
| 100 | 100 | 0.5000 | 0.5016 | ** | ** | 0.0001 | 0.0001 | ** | ** |

(*) The estimator is not available for $T = 2$.
(**) Computational cost is prohibitive for large $T$.

model). Most—but not all—conclusions drawn from Table 1 hold here. MILE continues to outperform the AB and AS estimators in terms of mean and MSE. If $\rho^* = -0.5$, MILE and BCOLS seem to perform similarly. If $\rho^* = 1.0$, MILE again performs better than BCOLS for small values of $T$.

**7. Conclusion.** A standard method to estimate parameters is the maximum likelihood estimator (MLE). In the presence of nuisance parameters, this approach concentrates out the likelihood by replacing these parameters with maximum likelihood estimators. An alternative approach entails maximizing a likelihood that depends only on parameters of interest. This marginal likelihood approach (e.g., [18] and [20]) yields an estimator for the structural parameter that is often less biased and more accurate than MLE (e.g., [11] and [24]).

TABLE 4

*Performance of estimators for the autoregressive parameter ρ (random effects, normal errors, and ρ = −0.50)*

| | | Mean | | | | MSE | | | |
|---|---|---|---|---|---|---|---|---|---|
| $T$ | $N$ | MILE | BCOLS | AB | AS | MILE | BCOLS | AB | AS |
| 2 | 5 | −0.5489 | −0.5689 | * | * | 0.1706 | 0.2478 | * | * |
| 2 | 10 | −0.5206 | −0.5622 | * | * | 0.0694 | 0.1020 | * | * |
| 2 | 25 | −0.5024 | −0.5485 | * | * | 0.0269 | 0.0374 | * | * |
| 2 | 100 | −0.5047 | −0.5476 | * | * | 0.0058 | 0.0104 | * | * |
| 3 | 5 | −0.4920 | −0.4907 | −0.0209 | −0.3722 | 0.0801 | 0.0791 | 20.5152 | 0.3044 |
| 3 | 10 | −0.5006 | −0.4994 | −0.4555 | −0.4485 | 0.0326 | 0.0352 | 4.0370 | 0.1651 |
| 3 | 25 | −0.5024 | −0.5087 | −0.4951 | −0.4990 | 0.0117 | 0.0146 | 0.0409 | 0.0578 |
| 3 | 100 | −0.5020 | −0.5063 | −0.4948 | −0.5368 | 0.0031 | 0.0033 | 0.0080 | 0.0129 |
| 5 | 5 | −0.4878 | −0.4728 | −0.5408 | −0.3755 | 0.0339 | 0.0371 | 0.0549 | 0.1201 |
| 5 | 10 | −0.4971 | −0.4871 | −0.5262 | −0.4113 | 0.0156 | 0.0202 | 0.0326 | 0.0713 |
| 5 | 25 | −0.5000 | −0.5007 | −0.5153 | −0.4608 | 0.0069 | 0.0073 | 0.0136 | 0.0310 |
| 5 | 100 | −0.4992 | −0.5021 | −0.5030 | −0.4860 | 0.0017 | 0.0017 | 0.0033 | 0.0069 |
| 10 | 5 | −0.4947 | −0.4779 | 0.6536 | −0.4602 | 0.0157 | 0.0181 | 3313.3070 | 0.0343 |
| 10 | 10 | −0.4965 | −0.4944 | −0.5334 | −0.4563 | 0.0083 | 0.0078 | 0.0098 | 0.0211 |
| 10 | 25 | −0.4987 | −0.4951 | −0.5144 | −0.4541 | 0.0031 | 0.0032 | 0.0046 | 0.0122 |
| 10 | 100 | −0.4995 | −0.4984 | −0.5024 | −0.4552 | 0.0008 | 0.0008 | 0.0014 | 0.0041 |
| 25 | 5 | −0.4958 | −0.4921 | ** | ** | 0.0061 | 0.0066 | ** | ** |
| 25 | 10 | −0.4986 | −0.4952 | ** | ** | 0.0033 | 0.0030 | ** | ** |
| 25 | 25 | −0.4988 | −0.4994 | ** | ** | 0.0013 | 0.0012 | ** | ** |
| 25 | 100 | −0.4996 | −0.4998 | ** | ** | 0.0003 | 0.0003 | ** | ** |
| 100 | 5 | −0.4996 | −0.4986 | ** | ** | 0.0016 | 0.0015 | ** | ** |
| 100 | 10 | −0.5002 | −0.4992 | ** | ** | 0.0008 | 0.0008 | ** | ** |
| 100 | 25 | −0.4997 | −0.4999 | ** | ** | 0.0003 | 0.0003 | ** | ** |
| 100 | 100 | −0.5000 | −0.4993 | ** | ** | 0.0001 | 0.0001 | ** | ** |

(*) The estimator is not available for $T = 2$.
(**) Computational cost is prohibitive for large $T$.

If the number of nuisance parameters increases, MLE may not even be consistent. This paper proposes a marginal likelihood approach to solve the incidental parameter problem. The use of invariance suggests which marginal likelihoods are to be maximized. We do not necessarily seek complete elimination of the incidental parameters. The goal is to find a group of transformations that preserves the structural parameters and yields a reduction in the incidental parameter space to a finite dimension.

We illustrate this approach with four examples: a stationary autoregressive model with fixed effects; a monotonic transformation model; an instrumental variable (IV) model; and a dynamic panel data model. In the first two examples, the invariant likelihoods are the products of marginal likelihoods and do not depend on the incidental parameters at all. In the last two examples, the invariant likelihoods are Wishart and depend on the incidental parameters through one-dimensional noncentrality parameters.

TABLE 5
*Performance of estimators for the autoregressive parameter $\rho$ (random effects, normal errors, and $\rho = 1.00$)*

| | | Mean | | | | MSE | | | |
|---|---|---|---|---|---|---|---|---|---|
| $T$ | $N$ | MILE | BCOLS | AB | AS | MILE | BCOLS | AB | AS |
| 2 | 5 | 0.9307 | 1.6990 | * | * | 0.1316 | 0.7595 | * | * |
| 2 | 10 | 0.9766 | 1.7115 | * | * | 0.0679 | 0.6034 | * | * |
| 2 | 25 | 1.0009 | 1.6943 | * | * | 0.0274 | 0.5166 | * | * |
| 2 | 100 | 0.9958 | 1.7047 | * | * | 0.0057 | 0.5048 | * | * |
| 3 | 5 | 0.9674 | 1.5029 | 1.0935 | 1.3267 | 0.0452 | 0.3211 | 36.9311 | 0.1953 |
| 3 | 10 | 1.0072 | 1.5032 | 1.0299 | 1.3320 | 0.0224 | 0.2776 | 5.5735 | 0.1386 |
| 3 | 25 | 0.9971 | 1.5156 | 1.0120 | 1.3469 | 0.0059 | 0.2733 | 0.0313 | 0.1318 |
| 3 | 100 | 0.9975 | 1.5216 | 0.9996 | 1.3624 | 0.0015 | 0.2740 | 0.0068 | 0.1345 |
| 5 | 5 | 0.9827 | 1.3241 | 0.9478 | 1.1497 | 0.0093 | 0.1190 | 0.0313 | 0.0363 |
| 5 | 10 | 0.9949 | 1.3341 | 0.9838 | 1.1531 | 0.0032 | 0.1165 | 0.0089 | 0.0289 |
| 5 | 25 | 0.9984 | 1.3403 | 0.9919 | 1.1659 | 0.0012 | 0.1174 | 0.0030 | 0.0294 |
| 5 | 100 | 0.9999 | 1.3442 | 0.9986 | 1.1760 | 0.0003 | 0.1189 | 0.0007 | 0.0315 |
| 10 | 5 | 0.9960 | 1.1774 | 1.2028 | 1.0534 | 0.0015 | 0.0330 | 55.2326 | 0.0065 |
| 10 | 10 | 0.9989 | 1.1838 | 0.9892 | 1.0621 | 0.0004 | 0.0343 | 0.0007 | 0.0053 |
| 10 | 25 | 0.9992 | 1.1839 | 0.9960 | 1.0680 | 0.0001 | 0.0340 | 0.0002 | 0.0051 |
| 10 | 100 | 1.0000 | 1.1854 | 0.9991 | 1.0687 | 0.0000 | 0.0344 | 0.0001 | 0.0048 |
| 25 | 5 | 0.9994 | 1.0765 | ** | ** | 0.0001 | 0.0059 | ** | ** |
| 25 | 10 | 1.0000 | 1.0767 | ** | ** | 0.0000 | 0.0059 | ** | ** |
| 25 | 25 | 0.9998 | 1.0776 | ** | ** | 0.0000 | 0.0060 | ** | ** |
| 25 | 100 | 1.0000 | 1.0776 | ** | ** | 0.0000 | 0.0060 | ** | ** |
| 100 | 5 | 1.0000 | 1.0197 | ** | ** | 0.0000 | 0.0004 | ** | ** |
| 100 | 10 | 0.9999 | 1.0198 | ** | ** | 0.0000 | 0.0004 | ** | ** |
| 100 | 25 | 1.0000 | 1.0198 | ** | ** | 0.0000 | 0.0004 | ** | ** |
| 100 | 100 | 1.0000 | 1.0198 | ** | ** | 0.0000 | 0.0004 | ** | ** |

(*) The estimator is not available for $T = 2$.
(**) Computational cost is prohibitive for large $T$.

For most groups of transformations, it is not possible to discard the incidental parameters completely. Because we allow invariant likelihoods to depend on incidental parameters, we have two considerations to make. First, finite-sample improvements may be possible using the orthogonalization approach of [15] to the invariant likelihood (e.g., [23]). Second, we treat the incidental parameters as an arbitrary sequence of numbers. Other authors (e.g., [21]) instead consider the incidental parameters as independently and identically distributed chance variables with distribution function. It would be interesting to understand the costs and benefits of treating the incidental parameters as unknown constants or chance variables.

## APPENDIX OF PROOFS

**Proofs of results stated in Section 3.**

PROOF OF LEMMA 3.1. Part (a) follows from Theorem 5.7 of [37]. Part (b) follows from Theorem 3.1 of [33]. Part (c) follows from Theorem 12.2.3 of [27] and Lemma 8.14 of [37].  □

PROOF OF PROPOSITION 3.1. For part (a), we need to show that $M(y_{i\cdot}) = M(\widetilde{y}_{i\cdot})$ if and only if $\widetilde{y}_{i\cdot} = y_{i\cdot} + \widetilde{g} \cdot 1_T$ for some $\widetilde{g}$. Clearly, $M(y_{i\cdot})$ is an invariant statistic,

$$M(y_{i\cdot} + g \cdot 1_T) = D(y_{i\cdot} + g \cdot 1_T) = Dy_{i\cdot} + g \cdot D1_T = Dy_{i\cdot} = M(y_{i\cdot}).$$

Now, suppose that $M(y_{i\cdot}) = M(\widetilde{y}_{i\cdot})$. This implies that $Dz_i = 0$ for $z_i = \widetilde{y}_{i\cdot} - y_{i\cdot}$, which means that $z_i$ belongs to the space orthogonal to the row space of $D$. Because $\mathrm{rank}(D) = T - 1$, the orthogonal space has dimension one. As this space contains the vector $1_T$, it must be the case that $z_i = \widetilde{g} \cdot 1_T$ for some scalar $\widetilde{g}$. Therefore, $\widetilde{y}_{i\cdot} = y_{i\cdot} + \widetilde{g} \cdot 1_T$.

Part (b) follows from the fact that the group of transformations acts transitively on $\eta_i$. Part (c) follows from the formula of the density of a normal distribution.  □

PROOF OF PROPOSITION 3.2. For part (a), let $M_{it}$ be the rank of $y_{it}$ in the collection $y_{i1}, \ldots, y_{iT}$. Formally, we can define $M_{it}$ through $y_{it} = y_{i(M_{it})}$. We shall abbreviate the notation, for example, $(g(y_{i1}), g(y_{i2}), \ldots, g(y_{iT}))$ as $g(y_{i\cdot})$. The maximal invariant is $M_i = (M_{i1}, \ldots, M_{iT}) = M(y_{i\cdot})$. We need to show that $M(y_{i\cdot}) = M(\widetilde{y}_{i\cdot})$ if and only if $\widetilde{y}_{i\cdot} = \widetilde{g}(y_{i\cdot})$. Consider the case that if $t \neq \widetilde{t}$, then $y_{it} \neq y_{i\widetilde{t}}$ (this set has probability measure equal to one). Clearly, $M_i$ is an invariant statistic. Now, suppose that $M(y_{i\cdot}) = M(\widetilde{y}_{i\cdot})$. This implies that $M_{i1} = \widetilde{M}_{i1}, \ldots, M_{iT} = \widetilde{M}_{iT}$. Therefore, $y_{ij_1} < \cdots < y_{ij_T}$ and $\widetilde{y}_{ij_1} < \cdots < \widetilde{y}_{ij_T}$. There is a continuous, strictly increasing transformation $\widetilde{g}$ such that $\widetilde{y}_{it} = \widetilde{g}(y_{it})$, $t = 1, \ldots, T$.

Part (b) follows from the fact that the group of transformations acts transitively on $\eta_i$.

For part (c), we note that because $\eta_i$ is an increasing transformation, $M_{it}$ is also the rank in the collection $y^*_{i1}, \ldots, y^*_{iT}$, where $y^*_{it} = x'_{it}\beta + u_{it}$. We note that $y^*_{i1}, \ldots, y^*_{iT}$ are jointly independent with marginal densities

$$f_{it}(z_{it}; \beta) = \frac{1}{\sqrt{2\pi}} \exp\left\{-\frac{1}{2}(z_{it} - x'_{it}\beta)^2\right\}.$$

Now, we note that

$$P(M_{i1} = m_{i1}, \ldots, M_{iT} = m_{iT})$$
$$= \int \cdots \int f_{i1}(z_{i1}; \beta) \cdots f_{iT}(z_{iT}; \beta) \, dz_{i1} \cdots dz_{iT},$$

integrated over the set in which $z_{it}$ is the $m_{it}$th smallest element of $z_{i1}, \ldots, z_{iT}$. We follow [27] and transform $w_{m_{it}} = z_{it}$ to obtain

$$P(M_{i1} = m_{i1}, \ldots, M_{iT} = m_{iT}) = \int_A \prod_{t=1}^T f_{it}(w_{m_{it}}; \beta) \, dw$$

$$= \int_A \prod_{t=1}^T \frac{f_{it}(w_{m_{it}}; \beta)}{f(w_{m_{it}})} f(w_{m_{it}}) \, dw,$$

where $f(w_t)$ is the density of a $N(0,1)$ distribution and $A = \{w \in \mathbb{R}^T; \; w_1 < \cdots < w_T\}$. Simple algebraic manipulations show that

$$P(M_i = m_i)$$

$$= \int_A \exp\left\{-\frac{1}{2}\sum_{t=1}^T (w_{m_{it}} - x'_{it}\beta)^2 + \frac{1}{2}\sum_{t=1}^T w_{m_{it}}^2\right\} \prod_{t=1}^T f(w_{m_{it}}) \, dw$$

$$= \int_A \exp\left\{\sum_{t=1}^T w_{m_{it}} x'_{it}\beta - \frac{1}{2}\sum_{t=1}^T (x'_{it}\beta)^2\right\} \prod_{t=1}^T f(w_{m_{it}}) \, dw$$

$$= \frac{1}{T!} \int_A \exp\left\{\left(\sum_{t=1}^T w_{m_{it}} x'_{it}\right)\beta - \frac{1}{2}\beta'\left(\sum_{t=1}^T x_{it} x'_{it}\right)\beta\right\} T! \prod_{t=1}^T f(w_{m_{it}}) \, dw,$$

where $T! \prod_{t=1}^T f(w_t)$ for $w_1 < \cdots < w_T$ is the p.d.f. of $V_{(1)}, \ldots, V_{(T)}$. $\quad\square$

**Proofs of results stated in Section 4.** For convenience, we omit the subscript in $\lambda_N$.

PROOF OF PROPOSITION 4.1. For part (a), we need to show that $M(R_1, R_2) = M(\widetilde{R}_1, \widetilde{R}_2)$, if and only if $(\widetilde{R}_1, \widetilde{R}_2) = (\widetilde{g}R_1, R_2)$ for some $\widetilde{g} \in O(K)$. Clearly, $M(y_i.)$ is an invariant statistic,

$$M(gR_1, R_2) = (R'_1 g' g R_1, R_2) = (R'_1 R_1, R_2) = M(R_1, R_2).$$

Now, suppose that $M(R_1, R_2) = M(\widetilde{R}_1, \widetilde{R}_2)$. This is equivalent to $R'_1 R_1 = \widetilde{R}'_1 \widetilde{R}_1$ and $R_2 = \widetilde{R}_2$. But this implies that $\widetilde{R}_1 = \widetilde{g}R_1$ (and, of course, $R_2 = \widetilde{R}_2$).

Part (b) follows analogously. $\quad\square$

PROOF OF THEOREM 4.1. Following [5], the density function of $Y' N_Z Y$ at $q$ is

$$f(q) = C_{1,K} \cdot \exp\left(-\frac{N\lambda}{2} a'\Sigma^{-1}a\right)|\Sigma|^{-K/2}|q|^{(K-3)/2}\exp\left(-\frac{1}{2}\operatorname{tr}(\Sigma^{-1}q)\right)$$

$$\times (\sqrt{N\lambda \cdot a'\Sigma^{-1}q\Sigma^{-1}a})^{-(K-2)/2} I_{(K-2)/2}(\sqrt{N\lambda \cdot a'\Sigma^{-1}q\Sigma^{-1}a}).$$

The density function of $W_N$ is then

$$g(w; \beta, \lambda_N) = f(q(w)) \cdot |q'(w)| = f(q(w))N^{2 \cdot 3/2},$$

which simplifies to (4.1). $\square$

PROOF OF THEOREM 4.2.    The log-likelihood function divided by $N$ is

$$Q_N(\theta) = -\frac{1}{2}\lambda \cdot a'\Sigma^{-1}a + \frac{1}{N}\ln\left(Z_N^{-(K-2)/2}I_{(K-2)/2}\left(\frac{N}{2}Z_N\right)\right)$$

$$(A.1) \qquad -\frac{K}{2N}\ln|\Sigma| + \frac{K-3}{2N}\ln|W_N| - \frac{1}{2}\operatorname{tr}(\Sigma^{-1}W_N)$$

$$+ \frac{1}{N}\ln(2^{(K-2)/2}N^{(K+2)/2}C_{1,K}),$$

where $Z_N = 2\sqrt{\lambda \cdot a'\Sigma^{-1}W_N\Sigma^{-1}a}$.

All terms in the last two lines converge under both SIV and MWIV asymptotics (the only exception is $\ln|W_N|$ under SIV asymptotics and under MWIV asymptotics with $\alpha = 0$). For example, the last term is

$$\frac{1}{N}\ln(2^{(K-2)/2}N^{(K+2)/2}C_{1,K}) = \frac{1}{N}\ln\left(\frac{N^{(K+2)/2}}{\Gamma((K-1)/2)}\right) + o(1)$$

under both SIV and MWIV asymptotics. Under SIV asymptotics,

$$\frac{1}{N}\ln\left(\frac{N^{(K+2)/2}}{\Gamma((K-1)/2)}\right) \to 0.$$

Under MWIV asymptotics, we can use Stirling's formula to obtain

$$\frac{1}{N}\ln\left(\frac{N^{(K+2)/2}}{\Gamma((K-1)/2)}\right) \to \frac{\alpha}{2}\left\{1 - \ln\left(\frac{\alpha}{2}\right)\right\}.$$

However, the second and third lines in (A.1) do not depend on $\theta$. As a result, these terms can be ignored in finding the limiting behavior of $\widehat{\theta}_N$. Hence, define the objective function

$$\widehat{Q}_N(\theta) = -\frac{1}{2}\lambda \cdot a'\Sigma^{-1}a + \frac{1}{N}\ln\left(Z_N^{-(K-2)/2}I_{(K-2)/2}\left(\frac{N}{2}Z_N\right)\right).$$

The quantity $Z_N$ depends on $W_N$. Following [32], Section 10.2,

$$E(W_N) = \frac{K \cdot \Sigma + \overline{M}'\overline{M}}{N} = \frac{K \cdot \Sigma + \pi'Z'Z\pi \cdot a^*a^{*\prime}}{N} = \frac{K}{N}\Sigma + \lambda_N^* \cdot a^*a^{*\prime}.$$

From here, we split the result into SIV or MWIV with $\alpha = 0$ asymptotics, and MWIV with $\alpha > 0$.

For part (a), $W_N = W_N^* + o_p(1)$, where

$$W_N^* \equiv \lambda_N^* \cdot a^*a^{*\prime}.$$

Hence, $Z_N = Z_N^* + o_p(1)$, where

$$Z_N^* \equiv 2\sqrt{\lambda \cdot \lambda_N^* (a'\Sigma^{-1}a^*)^2}.$$

The same holds for nonnormal errors, as long as $V(W_N) \to 0$.

Because $K$ is fixed and $N \to \infty$, $\widehat{Q}_N(\theta) = \overline{Q}_N(\theta) + o_p(1)$ (uniformly in $\theta \in \Theta$ compact), where

$$\overline{Q}_N(\theta) = -\frac{1}{2}\lambda \cdot a'\Sigma^{-1}a + \lambda^{1/2}\lambda_N^{*1/2}a^{*\prime}\Sigma^{-1}a.$$

The first-order condition (FOC) for $\overline{Q}_N(\theta)$ is given by

$$\frac{\partial \overline{Q}_N(\theta)}{\partial \beta} = -\lambda \cdot a'\Sigma^{-1}e_1 + \lambda^{1/2}\lambda_N^{*1/2}a^{*\prime}\Sigma^{-1}e_1,$$

$$\frac{\partial \overline{Q}_N(\theta)}{\partial \lambda} = -\frac{1}{2}a'\Sigma^{-1}a + \frac{1}{2}\lambda^{-1/2}\lambda_N^{*1/2}a^{*\prime}\Sigma^{-1}a.$$

The value $\theta^* = (\beta^*, \lambda_N^*)$ minimizes $\overline{Q}_N(\theta)$, setting the FOC to zero.

For parts (a)(i), (ii), $\overline{Q}_N(\theta) \to_p \overline{Q}(\theta)$, where

$$\overline{Q}(\theta) = -\frac{1}{2}\lambda \cdot a'\Sigma^{-1}a + \lambda^{1/2}\lambda^{*1/2}a^{*\prime}\Sigma^{-1}a.$$

Since $\theta \in \Theta$ compact and $\overline{Q}(\theta)$ is continuous, $\widehat{\theta}_N \to_p \theta$.

For part (a)(iii), we can define $\tau(\theta, \theta_N^*) \equiv \overline{Q}_N(\theta)$ which is continuous. For each point $\theta_N^*$, the function $\tau(\theta, \theta_N^*)$ reaches the maximum at $\theta = \theta_N^*$. Because $\theta \in \Theta$ compact and $\tau(\cdot, \theta_N^*)$ is continuous,

$$\sup_{\theta \in \Theta; \|\theta - \theta_N^*\| \geq \varepsilon} \overline{Q}_N(\theta) - \overline{Q}_N(\theta_N^*) = \max_{\theta \in \Theta; \|\theta - \theta_N^*\| \geq \varepsilon} \overline{Q}_N(\theta) - \overline{Q}_N(\theta_N^*) \equiv \delta(\theta_N^*) < 0.$$

Because $0 < \liminf \lambda_N^*$ and $\limsup \lambda_N^* < \infty$, there exists a compact set $\Theta^*$ such that $0 \notin \Theta^*$ in which $\theta_N^* \in \Theta^*$ eventually. Using continuity of $\delta(\cdot)$,

$$\sup_{\theta_N^* \in \Theta^*} \delta(\theta_N^*) = \max_{\theta_N^* \in \Theta^*} \delta(\theta_N^*) = \delta < 0$$

for large enough $N$. This implies $\theta_N^*$ is an identifiably unique sequence of maximizers of $\overline{Q}_N(\theta)$,

$$\limsup \sup_{\theta \in \Theta; \|\theta - \theta_N^*\| \geq \varepsilon} \overline{Q}_N(\theta) - \overline{Q}_N(\theta_N^*) < 0.$$

The result now follows from [36], Lemma 3.1.

For part (b), $W_N = W_N^* + o_p(1)$ under SIV and MWIV asymptotics, where

$$W_N^* = \alpha\Sigma + \lambda_N^* \cdot a^* a^{*\prime}.$$

Hence, $Z_N = Z_N^* + o_p(1)$, where $Z_N^*$ is defined as

$$Z_N^* \equiv 2\sqrt{\lambda \cdot a'\Sigma^{-1}(\alpha\Sigma + \lambda_N^* \cdot a^* a^{*\prime})\Sigma^{-1}a}.$$

The same holds for nonnormal errors, as long as $V(W_N) \to 0$. For $K/N \to \alpha > 0$, we use [1] to show that $\widehat{Q}_N(\theta) = \overline{Q}_N(\theta) + o_p(1)$ (uniformly in $\theta \in \Theta$ compact), where

$$\overline{Q}_N(\theta) = -\frac{1}{2}\lambda \cdot a'\Sigma^{-1}a + \frac{\alpha}{2}\left(1 + \frac{Z_N^{*2}}{\alpha^2}\right)^{1/2} - \frac{\alpha}{2}\ln\left(1 + \left(1 + \frac{Z_N^{*2}}{\alpha^2}\right)^{1/2}\right).$$

The first-order condition (FOC) for $\overline{Q}_N(\theta)$ is given by

$$\frac{\partial \overline{Q}_N(\theta)}{\partial \beta} = -\lambda \cdot a'\Sigma^{-1}e_1 + \frac{2\alpha}{\alpha}\frac{\alpha \cdot a'\Sigma^{-1}e_1 + \lambda_N^* \cdot a^{*\prime}\Sigma^{-1}a \cdot a^{*\prime}\Sigma^{-1}e_1}{1 + (1 + Z_N^{*2}/\alpha^2)^{1/2}},$$

$$\frac{\partial \overline{Q}_N(\theta)}{\partial \lambda} = -\frac{1}{2}a'\Sigma^{-1}a + \frac{1}{\alpha}\frac{\alpha \cdot a'\Sigma^{-1}a + \lambda_N^* \cdot (a^{*\prime}\Sigma^{-1}a)^2}{1 + (1 + Z_N^{*2}/\alpha^2)^{1/2}}.$$

The value $\theta_N^* = (\beta^*, \lambda_N^*)$ minimizes $\overline{Q}_N(\theta)$, setting the FOC to zero.

For parts (b)(i), (ii), $\overline{Q}_N(\theta) \to_p \overline{Q}(\theta)$ given by

$$\overline{Q}(\theta) = -\frac{1}{2}\lambda \cdot a'\Sigma^{-1}a + \frac{\alpha}{2}\left(1 + \frac{Z_N^{*2}}{\alpha^2}\right)^{1/2} - \frac{\alpha}{2}\ln\left(1 + \left(1 + \frac{Z_N^{*2}}{\alpha^2}\right)^{1/2}\right),$$

where $Z^* \equiv 2\sqrt{\lambda \cdot a'\Sigma^{-1}(\alpha\Sigma + \lambda^* \cdot a^* a^{*\prime})\Sigma^{-1}a}$. Since $\theta \in \Theta$ compact and $\overline{Q}(\theta)$ is continuous, $\widehat{\theta}_N \to_p \theta$.

Part (b)(iii) follows analogously to part (a)(iii).  $\square$

PROOF OF PROPOSITION 4.2.  It follows from [12] that the integrated likelihood [over Haar measures for $O(k)$] is maximized over $a$ by

$$\max_a \frac{a'\Sigma^{-1/2}Y'N_Z Y\Sigma^{-1/2}a}{a'a}.$$

This optimal $a$ is the eigenvector corresponding to the largest eigenvalue of $\Sigma^{-1/2}Y'N_Z Y\Sigma^{-1/2}$. The integrated likelihood coincides with the likelihood of the maximal invariant and $a$ is a transformation of $\beta$. As a result, MILE is equivalent to LIMLK.  $\square$

PROOF OF THEOREM 4.3.  For part (a), when $K$ is fixed or $K/N \to 0$,

(A.2)   $\widehat{Q}_N(\theta) = -\frac{1}{2}\lambda \cdot a'\Sigma^{-1}a + \lambda^{1/2}(a'\Sigma^{-1}W_N\Sigma^{-1}a)^{1/2} + o_p(N^{-1}).$

All results below hold up to $o_p(N^{-1/2})$ order.

The components of the score function $S_N(\theta)$ are

$$\frac{\partial Q_N(\theta)}{\partial \beta} = -\lambda \cdot a'\Sigma^{-1}e_1 + \lambda^{1/2}\frac{a'\Sigma^{-1}W_N\Sigma^{-1}e_1}{(a'\Sigma^{-1}W_N\Sigma^{-1}a)^{1/2}},$$

$$\frac{\partial Q_N(\theta)}{\partial \lambda} = -\frac{a'\Sigma^{-1}a}{2} + \frac{(a'\Sigma^{-1}W_N\Sigma^{-1}a)^{1/2}}{2\lambda^{1/2}}.$$

The components of the Hessian matrix $H_N(\theta) \equiv H(W_N;\theta)$ are

$$\frac{\partial^2 Q_N(\theta)}{\partial \beta^2} = -\lambda \cdot e_1'\Sigma^{-1}e_1 + \lambda^{1/2}\frac{e_1'\Sigma^{-1}W_N\Sigma^{-1}e_1}{(a'\Sigma^{-1}W_N\Sigma^{-1}a)^{1/2}}$$

$$- \lambda^{1/2}\frac{(a'\Sigma^{-1}W_N\Sigma^{-1}e_1)^2}{(a'\Sigma^{-1}W_N\Sigma^{-1}a)^{3/2}},$$

$$\frac{\partial^2 Q_N(\theta)}{\partial \beta \, \partial \lambda} = -a'\Sigma^{-1}e_1 + \frac{a'\Sigma^{-1}W_N\Sigma^{-1}e_1}{2\lambda^{1/2}(a'\Sigma^{-1}W_N\Sigma^{-1}a)^{1/2}},$$

$$\frac{\partial^2 Q_N(\theta)}{\partial \lambda^2} = -\frac{1}{4}\frac{(a'\Sigma^{-1}W_N\Sigma^{-1}e_1)^{1/2}}{\lambda^{3/2}}.$$

Because $W_N \to_p W^*$, $H_N(\theta) \to_p -\mathcal{I}_0(\theta^*)$. Furthermore, $H_N(\theta) \to_p H(W_N^*;\theta)$ uniformly on $\theta = (\beta, \lambda)$ for a compact set containing $\theta^*$, as long as $\lambda > 0$. This completes part (a)(ii). To show part (a)(i), we write

$$\sqrt{N}S_N(\theta^*) \equiv \sqrt{N}S(W_N;\theta^*) \equiv \sqrt{N}[S(W_N;\theta^*) - S(W^*;\theta^*)].$$

Using $\mathrm{vec}(W_N) = \mathcal{D}_T\,\mathrm{vech}(W_N)$, where $\mathcal{D}_T$ is the duplication matrix (e.g., [30]), we write

$$\sqrt{N}S_N(\theta^*) \equiv \sqrt{N}[L(\mathrm{vech}(W_N);\theta^*) - L(\mathrm{vech}(W^*);\theta^*)],$$

where $L: \mathbf{R}^3 \to \mathbf{R}^2$. Now, $\sqrt{N}(\mathrm{vech}(W_N) - \mathrm{vech}(W^*))$ converges to a normal distribution by a standard CLT. As a result, using the delta method and the information identity, $\sqrt{N}S_N(\theta^*)$ converges to a normal distribution with zero mean and variance $\mathcal{I}_0(\theta^*)$. Part (iii) follows from [33].

For part (b), when $K/N \to \alpha > 0$,

$$(A.3) \quad \widehat{Q}_N(\theta) = -\frac{1}{2}\lambda \cdot a'\Sigma^{-1}a + \frac{\alpha}{2}\left(1 + \frac{Z_N^2}{\alpha^2}\right)^{1/2} - \frac{\alpha}{2}\ln\left(1 + \left(1 + \frac{Z_N^2}{\alpha^2}\right)^{1/2}\right)$$

up to an $o_p(N^{-1})$ term. All results below hold up to $o_p(N^{-1/2})$ order.

The components of the score function $S_N(\theta)$ are

$$\frac{\partial Q_N(\theta)}{\partial \beta} = -\lambda \cdot a'\Sigma^{-1}e_1 + \frac{2\lambda}{\alpha}\frac{a'\Sigma^{-1}W_N\Sigma^{-1}e_1}{1 + (1 + Z_N^2/\alpha^2)^{1/2}},$$

$$\frac{\partial Q_N(\theta)}{\partial \lambda} = -\frac{a'\Sigma^{-1}a}{2} + \frac{1}{\alpha}\frac{a'\Sigma^{-1}W_N\Sigma^{-1}a}{1 + (1 + Z_N^2/\alpha^2)^{1/2}}.$$

The components of the Hessian matrix $H_N(\theta)$ are

$$\frac{\partial^2 Q_N(\theta)}{\partial \beta^2} = -\lambda \cdot e_1' \Sigma^{-1} e_1 + \frac{2\lambda}{\alpha} \frac{e_1' \Sigma^{-1} W_N \Sigma^{-1} e_1}{1 + (1 + Z_N^2/\alpha^2)^{1/2}}$$

$$- \frac{8\lambda^2}{\alpha^3 (1 + Z_N^2/\alpha^2)^{1/2}} \frac{(a' \Sigma^{-1} W_N \Sigma^{-1} e_1)^2}{(1 + (1 + Z_N^2/\alpha^2)^{1/2})^2},$$

$$\frac{\partial^2 Q_N(\theta)}{\partial \beta \, \partial \lambda} = -a' \Sigma^{-1} e_1 + \frac{2}{\alpha} \frac{a' \Sigma^{-1} W_N \Sigma^{-1} e_1}{1 + (1 + Z_N^2/\alpha^2)^{1/2}}$$

$$- \frac{4\lambda \cdot a' \Sigma^{-1} W_N \Sigma^{-1} e_1}{\alpha^3 (1 + Z_N^2/\alpha^2)^{1/2}} \frac{a' \Sigma^{-1} W_N \Sigma^{-1} a}{(1 + (1 + Z_N^2/\alpha^2)^{1/2})^2},$$

$$\frac{\partial^2 Q_N(\theta)}{\partial \lambda^2} = \frac{-2}{\alpha^3 (1 + Z_N^2/\alpha^2)^{1/2}} \frac{(a' \Sigma^{-1} W_N \Sigma^{-1} a)^2}{(1 + (1 + Z_N^2/\alpha^2)^{1/2})^2}.$$

Parts (b)(i)–(iii) follow analogously to parts (a)(i)–(iii). $\square$

Proof of Corollary 4.1. The determinant of $\mathcal{I}_\alpha(\theta^*)$ simplifies to

$$|\mathcal{I}_\alpha(\theta^*)| = \frac{\lambda^{*2}(a^{*\prime} \Sigma^{-1} a^*)^2}{\alpha + 2\lambda^* \cdot a^{*\prime} \Sigma^{-1} a^*} \frac{a^{*\prime} \Sigma^{-1} a^* \cdot e_1' \Sigma^{-1} e_1 - (a^{*\prime} \Sigma^{-1} e_1)^2}{2(\alpha + \lambda^* \cdot a^{*\prime} \Sigma^{-1} a^*)}.$$

Hence, the entry $(1, 1)$ of the inverse of $\mathcal{I}_\alpha(\theta^*)$ equals

$$(\mathcal{I}_\alpha(\theta^*)^{-1})_{11} = \frac{(a^{*\prime} \Sigma^{-1} a^*)^2}{2(\alpha + 2\lambda^* a^{*\prime} \Sigma^{-1} a^*)} |\mathcal{I}_\alpha(\theta^*)|^{-1}$$

$$= \frac{\alpha + \lambda^* \cdot a^{*\prime} \Sigma^{-1} a^*}{\lambda^{*2} \cdot a^{*\prime} \Sigma^{-1} a^*} \frac{a^{*\prime} \Sigma^{-1} a^*}{a^{*\prime} \Sigma^{-1} a^* \cdot e_1' \Sigma^{-1} e_1' - (a^{*\prime} \Sigma^{-1} e_1)^2}$$

$$= \frac{\sigma_u^2}{\lambda^{*2}} \left\{ \lambda^* + \frac{\alpha}{a^{*\prime} \Sigma^{-1} a^*} \right\}.$$

This expression coincides with the asymptotic variance of LIMLK as described in (4.7) of [10]:

$$(\mathcal{I}_\alpha(\theta^*)^{-1})_{11} = \frac{\sigma_u^2}{\lambda^{*2}} \left\{ \lambda^* + \alpha \cdot e_2' \Sigma e_2 - \alpha \frac{(b' \Sigma e_2)^2}{b' \Sigma b} \right\}. \qquad \square$$

Proof of Theorem 4.4. This result follows from standard limit of experiment arguments (see [14]). Part (a) follows from expansions based on (A.2). Part (b) follows from expansions based on (A.3). $\square$

**Proofs of results stated in Section 5.** For convenience, we omit the subscript in $\lambda_N$. For the next proofs, define the following four quantities:

$$c_1 = \text{tr}(DB^* B^{*\prime} D') + \lambda_N^* 1_T' B^{*\prime} D' DB^* 1_T,$$

$$c_2 = 1'_T DB^* B^{*\prime} D' 1_T + \lambda_N^* (1'_T DB^* 1_T)^2,$$

$$c_3 = 1'_T F 1_T + (\rho^* - \rho) 1'_T F' F 1_T + \lambda^* 1'_T DB^* 1_T \cdot 1'_T F 1_T,$$

$$c_4 = (\rho^* - \rho) \operatorname{tr}(F'F) + \lambda^* \{ 1'_T F 1_T + (\rho^* - \rho) 1'_T F' F 1_T \}.$$

PROOF OF PROPOSITION 5.1. We omit the proof here as it has been generalized by [13]. □

PROOF OF THEOREM 5.1. The density function of $M$ at $q$ is

$$f(q) = C_{2,N} \cdot \exp\left(-\frac{\eta'\eta}{2\sigma^2} T\right) (\sigma^2)^{-NT/2} |q|^{(N-T-1)/2} \exp\left(-\frac{1}{2\sigma^2} \operatorname{tr}(DqD')\right)$$

$$\times \left(\sqrt{\frac{\eta'\eta}{(\sigma^2)^2} 1'_T DqD' 1_T}\right)^{-(N-2)/2} I_{(N-2)/2}\left(\sqrt{\frac{\eta'\eta}{(\sigma^2)^2} 1'_T DqD' 1_T}\right).$$

The density function of $W_N$ is then

$$g(w; \beta, \lambda_N) = f(q(w)) \cdot |q'(w)| = f(q(w)) N^{T(T+1)/2},$$

which simplifies to (5.3). □

PROOF OF THEOREM 5.2. The log-likelihood divided by $NT$ is

$$\begin{aligned}
Q_N(\theta) = &-\frac{1}{2} \ln \sigma^2 - \frac{1}{2\sigma^2} \frac{\operatorname{tr}(DW_N D')}{T} - \frac{1}{2}\lambda \\
(A.4) \qquad &+ \frac{1}{NT} \ln\left(Z_N^{-(N-2)/2} I_{(N-2)/2}\left(\frac{N}{2} Z_N\right)\right) \\
&+ \frac{N-T-1}{2NT} \ln |W_N| + \frac{1}{NT} \ln(2^{(N-2)/2} N^{NT/2-(N-2)/2} C_{2,N}),
\end{aligned}$$

where $Z_N = 2\sqrt{\lambda \frac{1'_T DW_N D' 1_T}{\sigma^2}}$.

The third line is well-behaved when $N \to \infty$ with $T$ fixed. For example, using Stirling's formula,

$$\begin{aligned}
&\frac{1}{NT} \ln(2^{(N-2)/2} N^{NT/2-(N-2)/2} C_{2,N}) \\
&= \frac{1}{T} \ln\left(\frac{N^{NT/2-(N-2)/2} 2^{1/2}}{\prod_{t=1}^{T-1}(N-t)^{(N-t-1)/(2N)} \exp(-(N-t)/(2N))}\right) + o(1) \\
&= \frac{\ln(2)}{2T} - \frac{1}{T} \ln\left(\prod_{t=1}^{T-1}\left(1 - \frac{t}{N}\right)^{1/2} \exp\left(-\frac{1}{2}\right)\right) + o(1) \\
&= \frac{\ln(2)}{2T} + \frac{T-1}{2T} + o(1).
\end{aligned}$$

In addition, $W_N = W_N^* + o_p(1)$, where

$$W_N^* \equiv \sigma^{*2} B^* (I_T + \lambda_N^* 1_T 1_T') B^{*\prime} = \frac{N \cdot \Sigma + \overline{M'} \overline{M}}{N} = E(W_N).$$

Now,

$$|W_N^*| = |B^*| \cdot |\sigma^{*2}(I_T + \lambda_N^* 1_T 1_T')| \cdot |B^{*\prime}| = (\sigma^{*2})^T |I_T + \lambda_N^* 1_T 1_T'|$$
$$= (\sigma^{*2})^T (1 + \lambda_N^* T).$$

As a result, $\ln(W_N) = T \ln(\sigma^{*2}) + \ln(1 + \lambda_N^* T) + o_p(1)$.

It is unknown whether the third line in (A.4) is well-behaved with $T \to \infty$. However, since it does not depend on $\theta$, it can be ignored when finding the limiting behavior of $\widehat{\theta}_N$. Hence, define the objective function

$$\widehat{Q}_N(\theta) = -\frac{1}{2} \ln \sigma^2 - \frac{1}{2\sigma^2} \frac{\operatorname{tr}(D W_N D')}{T} - \frac{1}{2} \lambda$$
$$+ \frac{1}{NT} \ln \left( Z_N^{-(N-2)/2} I_{(N-2)/2} \left( \frac{N}{2} Z_N \right) \right).$$

From here, we split the result into fixed $T$ and large $T$ asymptotics. For part (a), in which $N \to \infty$ with $T$ fixed, $Z_N = Z_N^* + o_p(1)$, where

$$Z_N^* \equiv 2 \sqrt{\lambda \frac{1_T' D W_N^* D' 1_T}{\sigma^2}}.$$

We use [1] to show that $\widehat{Q}_N(\theta) = \overline{Q}_N(\theta) + o_p(1)$, where

$$\overline{Q}_N(\theta) = -\frac{1}{2} \ln \sigma^2 - \frac{1}{2\sigma^2} \frac{\operatorname{tr}(D W_N^* D')}{T} - \frac{1}{2} \lambda + \frac{1}{2T} (1 + Z_N^{*2})^{1/2}$$
$$- \frac{1}{2T} \ln(1 + (1 + Z_N^{*2})^{1/2}).$$

The first-order condition (FOC) for $\overline{Q}_N(\theta)$ is given by

$$\frac{\partial \overline{Q}_N(\theta)}{\partial \rho} = \frac{\sigma^{*2}}{\sigma^2} \frac{(\rho^* - \rho) \operatorname{tr}(FF') + \lambda^* \{1_T' F 1_T + (\rho^* - \rho) 1_T' F' F 1_T\}}{T}$$
$$- \frac{2\sigma^{*2}}{\sigma^2} \frac{\lambda}{1 + (1 + Z_N^{*2})^{1/2}}$$
$$\times \frac{1_T' F 1_T + (\rho^* - \rho) 1_T' F' F 1_T + \lambda^* (T + (\rho^* - \rho) 1_T' F 1_T) 1_T' F 1_T}{T},$$

$$\frac{\partial \overline{Q}_N(\theta)}{\partial \sigma^2} = -\frac{1}{2\sigma^2} + \frac{\sigma^{*2}}{2(\sigma^2)^2} \frac{c_1}{T} - \frac{\sigma^{*2}}{(\sigma^2)^2} \frac{\lambda_N^*}{1 + (1 + Z_N^{*2})^{1/2}} \frac{c_2}{T},$$

$$\frac{\partial \overline{Q}_N(\theta)}{\partial \lambda} = -\frac{1}{2} + \frac{\sigma^{*2}}{\sigma^2} \frac{1}{1 + (1 + Z_N^{*2})^{1/2}} \frac{c_2}{T}.$$

The value $\theta^* = (\rho^*, \sigma^{*2}, \lambda_N^*)$ minimizes $\overline{Q}_N(\theta)$, setting the FOC to zero.

For parts (a)(i), (ii), $\overline{Q}_N(\theta) \to_p \overline{Q}(\theta)$ (uniformly in $\Theta$ compact) given by

$$\overline{Q}(\theta) = -\frac{1}{2}\ln\sigma^2 - \frac{1}{2\sigma^2}\frac{\text{tr}(DW^*D')}{T} - \frac{1}{2}\lambda + \frac{1}{2T}(1+Z^{*2})^{1/2}$$
$$- \frac{1}{2T}\ln(1+(1+Z^{*2})^{1/2}),$$

where $W^*$ and $Z^*$ are defined as

$$(A.5) \quad W^* = \sigma^{*2}B^*(I_T + \lambda^* 1_T 1_T')B^{*\prime} \quad \text{and} \quad Z^* = 2\sqrt{\lambda\frac{1_T'DW^*D'1_T}{\sigma^2}}.$$

Since $\theta \in \Theta$ compact and $\overline{Q}(\theta)$ is continuous, $\widehat{\theta}_N \to_p \theta$.

Part (a)(iii) follows analogously to Theorem 4.2(a)(iii).

For part (b), the dimension of $W_N$ changes as $T \to \infty$. Yet, for $|\rho^*| < 1$,

$$\frac{\text{tr}(DW_ND')}{T} = \lim_{T\to\infty}\frac{\text{tr}(DW_N^*D')}{T} + o_p(1)$$

and

$$\frac{1_T'DW_ND'1_T}{T^2} = \lim_{T\to\infty}\frac{1_T'DW_N^*D'1_T}{T^2} + o_p(1).$$

This approximation does not depend on how $N$ grows with $T$. We use [1] to obtain $\widehat{Q}_N(\theta) = \overline{Q}_N(\theta) + o_p(1)$, where

$$\overline{Q}_N(\theta) = -\frac{1}{2}\ln\sigma^2 - \frac{1}{2\sigma^2}\lim_{T\to\infty}\frac{\text{tr}(DW_N^*D')}{T} - \frac{1}{2}\lambda + \frac{1}{2}\lim_{T\to\infty}\frac{Z_N^*}{T}.$$

The first-order condition (FOC) for $\overline{Q}_N(\theta)$ is given by

$$\frac{\partial\overline{Q}_N(\theta)}{\partial\rho} = \lim_{T\to\infty}\frac{\sigma^{*2}}{\sigma^2}\frac{(\rho^*-\rho)\,\text{tr}(FF') + \lambda^*\{1_T'F1_T + (\rho^*-\rho)1_T'F'F1_T\}}{T}$$
$$- \lim_{T\to\infty}\frac{(\sigma^{*2})^{1/2}\lambda^{*1/2}\lambda^{1/2}}{(\sigma^2)^{1/2}}\frac{1_T'F1_T}{T},$$

$$\frac{\partial\overline{Q}_N(\theta)}{\partial\sigma^2} = -\frac{1}{2\sigma^2} + \lim_{T\to\infty}\frac{\sigma^{*2}}{2(\sigma^2)^2}\frac{c_1}{T} - \lim_{T\to\infty}\frac{(\sigma^{*2})^{1/2}\lambda^{1/2}\lambda^{*1/2}}{2(\sigma^2)^{3/2}}\frac{1_T'DB^*1_T}{T},$$

$$\frac{\partial\overline{Q}_N(\theta)}{\partial\lambda} = -\frac{1}{2} + \lim_{T\to\infty}\frac{(\sigma^{*2})^{1/2}\lambda^{*1/2}}{2(\sigma^2)^{1/2}\lambda^{1/2}}\frac{1_T'DB^*1_T}{T}.$$

The value $\theta^* = (\rho^*, \sigma^{*2}, \lambda_N^*)$ minimizes $\overline{Q}_N(\theta)$, setting the FOC to zero.

For parts (b)(i), (ii), $\widehat{Q}_N(\theta) = \overline{Q}(\theta) + o_p(1)$ (uniformly in $\Theta$ compact), given by

$$\overline{Q}(\theta) = -\frac{1}{2}\ln\sigma^2 - \frac{1}{2\sigma^2}\lim_{T\to\infty}\frac{\text{tr}(DW^*D')}{T} - \frac{1}{2}\lambda + \frac{1}{2}\lim_{T\to\infty}\frac{Z^*}{T},$$

where $W^*$ and $Z^*$ are defined in (A.5). Since $\theta \in \Theta$ compact and $\overline{Q}(\theta)$ is continuous, $\widehat{\theta}_N \to_p \theta$.

Part (b)(iii) follows analogously to Theorem 4.2(a)(iii).   $\square$

PROOF OF THEOREM 5.3.   First, we prove part (a). The objective function is

(A.6)
$$\widehat{Q}_N(\theta) = -\frac{\ln \sigma^2}{2} - \frac{\operatorname{tr}(DW_N D')}{2\sigma^2 T} - \frac{\lambda}{2} + \frac{(1 + Z_N^2)^{1/2}}{2T}$$
$$- \frac{\ln(1 + (1 + Z_N^2)^{1/2})}{2T}$$

up to an $o_p(N^{-1})$ term. All results below hold up to $o_p(N^{-1/2})$ order.

The components of the score function $S_N(\theta)$ are

$$\frac{\partial Q_N(\theta)}{\partial \rho} = \frac{1}{\sigma^2} \frac{\operatorname{tr}(J_T W_N D')}{T} - \frac{2\lambda}{1 + (1 + Z_N^2)^{1/2}} \frac{1_T' J_T W_N D' 1_T}{T},$$

$$\frac{\partial Q_N(\theta)}{\partial \sigma^2} = -\frac{1}{2\sigma^2} + \frac{1}{2(\sigma^2)^2} \frac{\operatorname{tr}(DW_N D')}{T}$$
$$- \frac{1}{(\sigma^2)^2} \frac{\lambda}{1 + (1 + Z_N^2)^{1/2}} \frac{1_T' DW_N D' 1_T}{T},$$

$$\frac{\partial Q_N(\theta)}{\partial \lambda} = -\frac{1}{2} + \frac{1}{\sigma^2} \frac{1}{1 + (1 + Z_N^2)^{1/2}} \frac{1_T' DW_N D' 1_T}{T}.$$

The Hessian matrix $H_N(\theta) \to_p -\mathcal{I}_T(\theta)$, whose components are

$$\frac{\partial^2 \overline{Q}_N(\theta)}{\partial \rho^2} = \frac{\sigma^{*2}}{\sigma^2} \frac{2\lambda}{1 + (1 + Z_N^{*2})^{1/2}} \frac{1_T' F' F 1_T + \lambda (1_T' F 1_T)^2}{T}$$
$$- \frac{\sigma^{*2}}{\sigma^2} \frac{\operatorname{tr}(F'F) + \lambda^* 1_T' F' F 1_T}{T}$$
$$- \left(\frac{\sigma^{*2}}{\sigma^2}\right)^2 \frac{8\lambda^2}{(1 + (1 + Z_N^{*2})^{1/2})^2} \frac{1}{(1 + Z_N^{*2})^{1/2}} \frac{(c_3)^2}{T},$$

$$\frac{\partial_N^2 \overline{Q}(\theta)}{\partial \rho \, \partial \sigma^2} = -\frac{\sigma^{*2}}{(\sigma^2)^2} \frac{c_4}{T} + \frac{\sigma^{*2}}{(\sigma^2)^2} \frac{2\lambda}{1 + (1 + Z_N^{*2})^{1/2}} \frac{c_3}{T}$$
$$\times \left\{1 - \frac{\sigma^{*2}}{\sigma^2} \frac{2\lambda c_2}{1 + (1 + Z_N^{*2})^{1/2}} \frac{1}{(1 + Z_N^{*2})^{1/2}}\right\},$$

$$\frac{\partial_N^2 \overline{Q}(\theta)}{\partial \rho \, \partial \lambda} = -\frac{\sigma^{*2}}{\sigma^2} \frac{2}{1 + (1 + Z_N^{*2})^{1/2}} \frac{c_3}{T}$$

$$\times \left\{ 1 - \frac{\sigma^{*2}}{\sigma^2} \frac{2\lambda c_2}{1 + (1 + Z_N^{*2})^{1/2}} \frac{1}{(1 + Z_N^{*2})^{1/2}} \right\},$$

$$\frac{\partial_N^2 \overline{Q}(\theta)}{\partial (\sigma^2)^2} = -\frac{(\sigma^{*2})^2}{(\sigma^2)^4} \frac{2\lambda^2}{(1 + (1 + Z_N^{*2})^{1/2})^2} \frac{1}{(1 + Z_N^{*2})^{1/2}} \frac{(c_2)^2}{T}$$

$$+ \frac{1}{2(\sigma^2)^2} - \frac{\sigma^{*2}}{(\sigma^2)^3} \frac{c_1}{T} + \frac{\sigma^{*2}}{(\sigma^2)^3} \frac{2\lambda}{1 + (1 + Z_N^{*2})^{1/2}} \frac{c_2}{T},$$

$$\frac{\partial_N^2 \overline{Q}(\theta)}{\partial \sigma^2 \partial \lambda} = -\frac{\sigma^{*2}}{(\sigma^2)^2} \frac{1}{1 + (1 + Z_N^{*2})^{1/2}} \frac{c_2}{T} \left\{ 1 - \frac{\sigma^{*2}}{\sigma^2} \frac{2\lambda c_2}{1 + (1 + Z_N^{*2})^{1/2}} \frac{1}{(1 + Z_N^{*2})^{1/2}} \right\},$$

$$\frac{\partial_N^2 \overline{Q}(\theta)}{\partial \lambda^2} = -\left( \frac{\sigma^{*2}}{\sigma^2} \right)^2 \frac{2}{(1 + (1 + Z_N^{*2})^{1/2})^2} \frac{1}{(1 + Z_N^{*2})^{1/2}} \frac{(c_2)^2}{T}.$$

This convergence is uniform on $\theta = (\beta, \lambda)$ for a compact set containing $\theta^*$, as long as $\lambda > 0$. This completes part (a)(ii). To show part (a)(i), we write

$$\sqrt{NT} S_N(\theta^*) \equiv \sqrt{NT} S(W_N; \theta^*) \equiv \sqrt{NT} [S(W_N; \theta^*) - S(W^*; \theta^*)].$$

Using $\text{vec}(W_N) = \mathcal{D}_T \text{vech}(W_N)$, where $\mathcal{D}_T$ is the duplication matrix (e.g., [30]), we write

$$\sqrt{NT} S_N(\theta^*) \equiv \sqrt{NT} [L(\text{vech}(W_N); \theta^*) - L(\text{vech}(W^*); \theta^*)],$$

where $L: \mathbf{R}^{T(T+1)/2} \to \mathbf{R}^3$. Now, $\sqrt{NT}(\text{vech}(W_N) - \text{vech}(W^*))$ converges to a normal distribution by a standard CLT. As a result, using the delta method and the information identity, $\sqrt{NT} S_N(\theta^*)$ converges to a normal distribution with zero mean and variance $\mathcal{I}_T(\theta^*)$. Part (iii) follows from [33].

Part (b) follows from the asymptotic normality of the score (whose variance is given by the reciprocal of the inverse of the limit of the Hessian matrix). As the remainder terms from expansions based on (A.6) are asymptotically negligible, (5.4) holds true. $\square$

PROOF OF COROLLARY 5.1. As a preliminary result, we need to find the limits of $T^{-1} \text{tr}(FF')$, $T^{-1} 1_T' F 1_T$ and $T^{-1} 1_T' F' F 1_T$, as $T \to \infty$. For the first term,

$$\frac{1}{T} \text{tr}(FF') = \frac{1}{T} \sum_{j=0}^{T-2} \sum_{i=0}^{j} \rho^{*2i} = \frac{T-1}{T} \sum_{i=0}^{T-1} \rho^{*2i} - \frac{1}{T} \sum_{i=0}^{T-1} i \rho^{*2i} \to \frac{1}{1 - \rho^{*2}},$$

because $\sum_{i=0}^{T-1} i(\rho^{*2})^i$ is a convergent series. This is true because a sufficient condition for a series $\sum_{i=0}^{T} a_i$ to converge is that $\lim \sqrt[T]{|a_T|} < 1$ as $T \to$ $\infty$. Taking $a_i = i(\rho^{*2})^i$, $\lim \sqrt[T]{|a_T|} = \lim \sqrt[T]{|T(\rho^{*2})^T|} = \rho^{*2} \lim \sqrt[T]{T} = \rho^{*2} < 1$. Analogously,

$$\frac{1}{T}1'_T F1_T = \frac{1}{T}\sum_{j=0}^{T-2}\sum_{i=0}^{j}\rho^{*i} = \frac{T-1}{T}\sum_{i=0}^{T-1}\rho^{*i} - \frac{1}{T}\sum_{i=0}^{T-1}i\rho^{*i} \to \frac{1}{1-\rho^*},$$

because $\sum_{i=0}^{T-1}i\rho^{*i}$ also converges. Finally, by the Cauchy–Schwarz inequality,

$$\left(\frac{1}{T}1'_T F1_T\right)^2 \le \frac{1}{T}1'_T F'F1_T = \frac{1}{T}\sum_{j=0}^{T-2}\left(\sum_{i=0}^{j}\rho^{*i}\right)^2 \le \frac{T-1}{T}\left(\frac{1}{1-\rho^*}\right)^2.$$

Taking limits, we obtain

$$\frac{1}{(1-\rho^*)^2} \le \liminf \frac{1}{T}1'_T F'F1_T \le \limsup \frac{1}{T}1'_T F'F1_T \le \frac{1}{(1-\rho^*)^2}.$$

Hence, the limit of $T^{-1}1'_T F'F1_T$ exists and equals $(1-\rho^*)^{-2}$.

Therefore, the limiting information matrix $\mathcal{I}_\infty(\theta^*)$ simplifies to

$$\mathcal{I}_\infty(\theta^*) = \begin{bmatrix} \dfrac{1}{1-\rho^{*2}} + \dfrac{\lambda^*}{(1-\rho^*)^2} & \dfrac{\lambda^*}{2\sigma^{*2}(1-\rho^*)} & \dfrac{1}{2(1-\rho^*)} \\ \dfrac{\lambda^*}{2\sigma^{*2}(1-\rho^*)} & \dfrac{2+\lambda^*}{4(\sigma^{*2})^2} & \dfrac{1}{4\sigma^{*2}} \\ \dfrac{1}{2(1-\rho^*)} & \dfrac{1}{4\sigma^{*2}} & \dfrac{1}{4\lambda^*} \end{bmatrix}.$$

The entry $(1,1)$ of the inverse of $\mathcal{I}_\infty(\theta^*)$ is

$$(\mathcal{I}_\infty(\theta^*)^{-1})_{11} = 1 - \rho^{*2}. \qquad \square$$

PROOF OF THEOREM 5.4. When $T \to \infty$, the objective function is

$$\widehat{Q}_N(\theta) = -\frac{1}{2}\ln\sigma^2 - \frac{1}{2\sigma^2}\frac{\mathrm{tr}(DW_N D')}{T} - \frac{1}{2}\lambda - \frac{1}{2T}Z_N$$

up to an $o_p(N^{-1})$ term. All results below hold up to $o_p(N^{-1/2})$ order.

The components of the score function $S_N(\theta)$ are

$$\frac{\partial Q_N(\theta)}{\partial\rho} = \frac{1}{\sigma^2}\frac{\mathrm{tr}(J_T W_N D')}{T} - \frac{\lambda^{1/2}}{(\sigma^2)^{1/2}}\frac{1'_T J_T W_N D'1_T}{T(1'_T DW_N D'1_T)^{1/2}},$$

$$\frac{\partial Q_N(\theta)}{\partial\sigma^2} = -\frac{1}{2\sigma^2} + \frac{1}{2(\sigma^2)^2}\frac{\mathrm{tr}(DW_N D')}{T} - \frac{\lambda^{1/2}}{2(\sigma^2)^{3/2}}\frac{(1'_T DW_N D'1_T)^{1/2}}{T},$$

$$\frac{\partial Q_N(\theta)}{\partial\lambda} = -\frac{1}{2} + \frac{1}{2(\sigma^2)^{1/2}\lambda^{1/2}}\frac{(1'_T DW_N D'1_T)^{1/2}}{T}.$$

If $|\rho^*|$ is bounded away from one, as $T \to \infty$,

$$\frac{\operatorname{tr}(J_T W_N D')}{T} \to_p \lim \frac{\operatorname{tr}(J_T W_N^* D')}{T},$$

$$\frac{1_T' J_T W_N D' 1_T}{T^2} \to_p \lim \frac{1_T' J_T W_N^* D' 1_T}{T^2},$$

$$\frac{\operatorname{tr}(D W_N D')}{T} \to_p \lim \frac{\operatorname{tr}(D W_N^* D')}{T}$$

and

$$\frac{1_T' D W_N D' 1_T}{T^2} \to_p \lim \frac{1_T' D W_N^* D' 1_T}{T^2}.$$

As a result, the Hessian matrix $-H_N(\theta) \to_p \mathcal{I}_\infty(\theta)$, whose components are limits of

$$-\frac{\partial^2 Q_N(\theta)}{\partial \rho^2} = \frac{\sigma^{*2}}{\sigma^2} \frac{\operatorname{tr}(F'F) + \lambda^* 1_T' F' F 1_T}{T},$$

$$-\frac{\partial^2 Q_N(\theta)}{\partial \rho \, \partial \sigma^2} = \frac{\sigma^{*2}}{(\sigma^2)^2} \frac{c_4}{T} - \frac{\lambda^{1/2} \lambda^{*1/2} (\sigma^{*2})^{1/2}}{2(\sigma^2)^{3/2}} \frac{1_T' F 1_T}{T},$$

$$-\frac{\partial^2 Q_N(\theta)}{\partial \rho \, \partial \lambda} = \frac{(\sigma^{*2})^{1/2} \lambda^{*1/2}}{2(\sigma^2)^{1/2} \lambda^{3/2}} \frac{1_T' F 1_T}{T},$$

$$-\frac{\partial^2 Q_N(\theta)}{\partial (\sigma^2)^2} = \frac{\sigma^{*2}}{(\sigma^2)^3} \frac{c_1}{T} - \frac{3}{4} \frac{(\sigma^{*2})^{1/2} \lambda^{1/2} \lambda^{*1/2}}{(\sigma^2)^{5/2}} \frac{1_T' D B^* 1_T}{T} - \frac{1}{2(\sigma^2)^2},$$

$$-\frac{\partial^2 Q_N(\theta)}{\partial \sigma^2 \, \partial \lambda} = \frac{(\sigma^{*2})^{1/2} \lambda^{*1/2}}{4(\sigma^2)^{3/2} \lambda^{1/2}} \frac{1_T' D B^* 1_T}{T} \quad \text{and}$$

$$-\frac{\partial^2 Q_N(\theta)}{\partial \lambda^2} = \frac{(\sigma^{*2})^{1/2} \lambda^{*1/2}}{4(\sigma^2)^{1/2} \lambda^{3/2}} \frac{1_T' D B^* 1_T}{T}.$$

This convergence is uniform on $\theta = (\beta, \lambda)$ for a compact set containing $\theta^*$, as long as $|\rho^*|$ is bounded away from one. This completes part (ii). To show part (i), define

$$\mathcal{W}_N = \left( \frac{\operatorname{tr}(J_T W_N D^{*\prime})}{T} \quad \frac{1_T' J_T W_N D^{*\prime} 1_T}{T^2} \quad \frac{\operatorname{tr}(D^* W_N' D^{*\prime})}{T} \quad \frac{1_T' D^* W_N' D^{*\prime} 1_T}{T^2} \right)'$$

and

$$\mathcal{W}_N^* = \left( \frac{\operatorname{tr}(J_T W_N^* D^{*\prime})}{T} \quad \frac{1_T' J_T W_N^* D^{*\prime} 1_T}{T^2} \quad \frac{\operatorname{tr}(D^* W_N^* D^{*\prime})}{T} \quad \frac{1_T' D^* W_N^* D^{*\prime} 1_T}{T^2} \right)'$$

and write

$$\sqrt{NT} S_N(\theta^*) \equiv \sqrt{NT}[L(\mathcal{W}_N; \theta^*) - L(\mathcal{W}_N^*; \theta^*)],$$

where $L : \mathbf{R}^4 \to \mathbf{R}^3$. Now, $\sqrt{NT}(\mathcal{W}_N - \mathcal{W}_N^*)$ converges to a normal distribution by a standard CLT and the Cramér–Wold device. Using the delta method and the information identity, $\sqrt{NT}S_N(\theta^*)$ converges to a normal distribution with zero mean and variance $\mathcal{I}_\infty(\theta^*)$, as long as $N \geq T$. Part (iii) follows from [33]. □

**Acknowledgments.** The author thanks Gary Chamberlain for helpful conversations and a correction, Rustam Ibragimov for valuable suggestions, Tiemen Woutersen for early discussions on the topic, an Associate Editor and a referee for several comments, and Jose Miguel Torres and Christiam Gonzales for research assistance.

Columbia University and FGV/EPGE
1022 International Affairs Building
MC 3308
420 West 118th Street
New York, New York 10027
USA
E-mail: mjmoreira@columbia.edu